\documentclass[11pt]{article}

\usepackage{amssymb,amsmath,amsthm}
\usepackage{graphicx}
\usepackage{multicol}

\usepackage{mathtools}
\mathtoolsset{showonlyrefs}

\numberwithin{equation}{section}

\usepackage{multirow}

\usepackage{color}

\def\N{\mathbb{N}}
\def\R{\mathbb{R}}
\def\Z{\mathbb{Z}}
\def\i{\mathrm i}
\def\d{\mathrm d}
\def\e{\mathrm e}
\def\E{\mathbb{E}}
\def\P{\mathbb{P}}
\def\text{\mbox}
\def\vep{\varepsilon}
\def\1{{\bf 1}}

\newcommand {\nn}{\nonumber}

\theoremstyle{plain}
\newtheorem{lemma}{Lemma}[section]
\newtheorem{proposition}[lemma]{Proposition}
\newtheorem{theorem}[lemma]{Theorem}

\theoremstyle{definition}
\newtheorem{example}[lemma]{Example}

\newtheorem{definition}[lemma]{Definition}
\newtheorem{remark}[lemma]{Remark}

\begin{document}

%

\title{Power variations for  fractional type   
	 infinitely divisible random fields
}

\date{\today}

\author{
	Andreas Basse-O'Connor\thanks{Department of Mathematics, Aarhus University, Denmark, E-mail: basse@math.au.dk.}, \
	Vytaut\.e Pilipauskait\.e\thanks{Department of Mathematics, University of Luxembourg, E-mail:  vytaute.pilipauskaite@gmail.com.}, \
	Mark Podolskij\thanks{Department of Mathematics, University of Luxembourg, E-mail: mark.podolskij@uni.lu.
	}
}

\maketitle

\begin{abstract}

This paper presents new limit theorems for power variations of fractional type symmetric infinitely divisible random fields. More specifically, the random field $X = (X(\boldsymbol{t}))_{\boldsymbol{t} \in [0,1]^d}$ is defined as an integral of a kernel function $g$ with respect to  a symmetric infinitely divisible random measure $L$ and is observed on a grid with mesh size $n^{-1}$. As $n \to \infty$, the first order limits are obtained for power variation statistics constructed from rectangular increments of $X$.
The present work is mostly related to \cite{BHP19,basse2017}, who studied a similar problem in 
the case $d=1$. 
We will see, however, that the asymptotic theory in the random field setting is much richer compared to  \cite{BHP19,basse2017} as it contains new limits, which depend on the precise structure of the kernel $g$.  We will 
give some important examples including the L\'evy moving average field, the well-balanced symmetric linear fractional $\beta$-stable sheet, and 
the moving average fractional $\beta$-stable field, and discuss potential consequences for statistical inference.

\end{abstract}



\section{Introduction}	 \label{sec1}
\setcounter{equation}{0}
\renewcommand{\theequation}{\thesection.\arabic{equation}}

The last decades have witnessed an immense progress  in limit theory for power variations of stochastic processes. 
Power variation functionals and related statistics play a major role in the analysis of the fine structure of the underlying model, in stochastic integration theory and statistical
applications. Asymptotic theory for power variations of various classes of stochastic processes has received a great deal of attention 
in the probabilistic and statistical literature.
We refer e.g. to \cite{BGJPS,J,JP,PV} for limit theory for power variations of It\^o semimartingales, to \cite{BNCP09,BNCPW09,Coeu,gl89,nr09} 
for the asymptotic results in the framework of fractional Brownian motion and related processes, and to \cite{Rosenblatt, CTV11, Rosenblatt-1} for investigations of power variation of the Rosenblatt process. 

More recently, there appeared numerous studies on limit theorems for statistics of non-Gaussian infinitely divisible moving-average processes. 
Central limit theorems for low frequency statistics of  infinite-variance stable moving averages have been investigated in \cite{PT2,PTA}. During the past years high frequency statistics of stationary increments L\'evy driven moving averages have been discussed in \cite{BHP19,basse2017}. In 
\cite{basse2017} the authors showed a variety of first and second order asymptotic results for power variation statistics, which heavily depend on the behaviour of the kernel near $0$, the Blumenthal--Getoor index of the driving L\'evy process and the considered power $p$. Later on these findings have been extended to a more general  class of statistics and processes in \cite{BHP18, BHP19}. We remark that the aforementioned probabilistic results are of immense importance for statistical applications. Indeed, they have been applied in \cite{LP19,LP20,MOP20} to obtain complete parametric estimation of the linear  fractional stable models and related processes in low and high frequency settings. Earlier studies on similar estimation problems, which are mainly concerned with estimation of the self-similarity parameter, can be found in \cite{AH,DI17,PTA,SPT}.  
Studies of high frequency statistics for L\'evy driven random fields are much more scarce in the literature.  Functional limit theorems for generalised 
variations of the fractional Brownian sheet have been investigated in \cite{PR16}, while power variations for certain integrals with respect to Gaussian white noise have been studied in \cite{pakk2014}.  We remark however that both classes of models are driven by a Gaussian field and the considered techniques do not apply in the more general L\'evy setting.

The aim of this paper is to study power variation statistics built from rectangular increments of certain random fields driven by an 
infinitely divisible random measure without a Gaussian part. More precisely,   we consider an  $\R$-valued random field
$X = (X(\boldsymbol{t}))_{\boldsymbol{t} \in \R^d}$ defined as
\begin{equation}\label{def:X}
( X(\boldsymbol{t}))_{\boldsymbol{t} \in \R^d} = \Big(  \int_{\R^d} g(\boldsymbol{t},\boldsymbol{s}) L(\d \boldsymbol{s}) \Big)_{\boldsymbol{t} \in \R^d},
\end{equation}
where $g:\R^d \times \R^d \to \R$ is a deterministic kernel to be introduced in \eqref{def:g} and $L$ is an infinitely divisible random measure on $\R^d$. 
We will focus on determining the first order asymptotic theory for power variation statistics of the form 
\begin{align} 
V_n(p) &:= \sum_{\boldsymbol{i} \in \{0,\dots,n-1\}^d
} | \Delta_{1/n} X(\boldsymbol{i}/n) |^p,  \\
\label{def:Del}   \Delta_{1/n} X(\boldsymbol{i}/n)&:= \sum_{\boldsymbol{\varepsilon} \in \{0,1\}^d} (-1)^{d+\sum_{j=1}^d \varepsilon_j} X
\left ((i_1+\varepsilon_1)/n, \dots,(i_d+\varepsilon_d)/n\right),
\end{align} 
where $\boldsymbol{i}=(i_1,\ldots, i_d)$, 
$\Delta_{1/n} X(\boldsymbol{i}/n)$ are rectangular increments of $X$, and $p>0$. We will show that the type of convergence and the limit of $V_n(p)$ crucially depend on the L\'evy measure of $L$, the considered power $p>0$ and the behaviour of rectangular increments of $g$ near  $\boldsymbol{0}\in \R^d$. These results can be considered as a  extension of  \cite[Theorem 1.1]{basse2017}
to the framework of random fields. However, the picture turns out to be more complex than for processes
studied in  \cite[Theorem 1.1]{basse2017}. Indeed, we will show that different forms of \textit{local homogeneity} of the kernel $g$, which are summarised in Assumptions (H1) 
and (H2), lead to different asymptotic results, a phenomenon that does not appear in the case $d=1$. 
In particular, the limit types stated 
in Theorems \ref{thm2}(i) and (ii) do not have a one-dimensional counterpart.     We will discuss how our theoretical results apply to most popular L\'evy driven random fields including the moving average field, the well-balanced symmetric linear fractional $\beta$-stable sheet and the moving average fractional $\beta$-stable field among other models. Furthermore, we will present a short discussion on potential application of our theory 
to parameter identification and parameter estimation.  

This paper is organised as follows. Section \ref{sec2} presents the model setting and the necessary definitions. The main theoretical results and their applications are demonstrated in  Section \ref{sec3}. All major proofs are collected in  Section \ref{sec4}. Some technical statements can be found in the Appendix.

\section{The setting, notations and definitions} \label{sec2}
\setcounter{equation}{0}
\renewcommand{\theequation}{\thesection.\arabic{equation}}

\subsection{Notations} \label{sec2.1}

Throughout the paper we denote all multi-indexed quantities by bold letters.
For $\boldsymbol{x} = (x_1,\dots,x_d) \in \R^d$ and $\boldsymbol{y} = (y_1,\dots,y_d) \in \R^d$, we write $\boldsymbol{x} < \boldsymbol{y}$ if $x_i < y_i$, $i=1,\dots,d$; the relation $\boldsymbol{x} \leq \boldsymbol{y}$ is defined similarly. 
	We denote the rectangle $[x_1, y_1] \times \cdots \times [x_d, y_d]$ by $[\boldsymbol{x}, \boldsymbol{y}]$ for  $\boldsymbol{x} \le \boldsymbol{y}$. 
	For each real number $x\in \R$ let $\{x\}=x-\lfloor x \rfloor$ denote its fractional part, and write $\{ \boldsymbol{x} \} = (\{x_1\}, \dots, \{x_d\} )$ for the fractional part of $\boldsymbol{x}\in \R^d$ taken coordinate-wise.   
	We set $\| \boldsymbol{x} \| = (x_1^2+\dots+x_d^2)^{1/2}$. 
	We define the open ball of radius $r > 0$ centered at a point $\boldsymbol{x}_0 \in \R^d$ as
	$B_r(\boldsymbol{x}_0) := \{ \boldsymbol{x} \in \R^d : \| \boldsymbol{x} - \boldsymbol{x}_0 \| < r  \}$. We denote
	the complement of a set $B$ in $\R^d$ by $B^c := \R^d \setminus B$. 
	Furthermore, ${\cal B}_b(\R^d)$ denotes a collection of all bounded Borel measurable subsets of $\R^d$ and  $\lambda^d$
	denotes the Lebesgue measure on $\R^d$.
	Finally, $\partial^d g (\boldsymbol{s})$ denotes the partial derivative $\frac{\partial^d}{\partial s_1 \dots \partial s_d} g (\boldsymbol{s})$ of $g$ at $\boldsymbol{s} \in \R^d$ if it exists, and otherwise we set $\partial^d g (\boldsymbol{s})$ equal to $0$.

We write $\overset{\P}{\to}$, $\overset{L^1}{\to}$, $\overset{\textnormal d}{\to}$ for convergence in probability, mean, distribution  of a sequence of random variables. The notation  $\overset{\textnormal d}{=}$ stands for equality in distribution of random variables and $\overset{\rm fdd}{=}$ denotes
the equality of finite-dimensional distributions of stochastic processes. We
write $Y_n \overset{{\cal F}\textnormal{-d}}{\to} Y$ if a sequence $( Y_n )_{n \in \N}$ of random variables defined on the probability space $(\Omega, {\cal F}, \P)$ converges ${\cal F}$-stably in law to $Y$. That is, $Y$ is a random variable defined on the extension of $(\Omega,{\cal F},\P)$ such that for all ${\cal F}$-measurable random variables $Z$ the joint convergence in distribution $(Y_n,Z) \overset{\textnormal d}{\to} (Y,Z)$ holds. 
For a detailed treatment of stable convergence we refer to \cite{hausler2015}.

Finally,  $C$ stands for a generic positive finite constant whose precise value is unimportant and may change from line to line. By convention, summation and product over an empty set is $0$ and $1$, respectively.

\subsection{The model} \label{sec2.2}

We consider a random field $X = (X(\boldsymbol{t}))_{\boldsymbol{t} \in \R^d}$ defined in
\eqref{def:X} as an integral of a kernel $g$ with respect to an infinitely divisible random measure $L$. 
We recall that the collection $L = (L(B))_{B \in {\cal B}_b (\R^d)}$ is an  infinitely divisible random measure when 
\begin{itemize}
	\item[(i)] for every sequence $( B_i )_{i \in \N}$ of pairwise disjoint sets in ${\cal B}_b (\R^d)$, $(L(B_i))_{i\in \N}$ forms a sequence of independent random variables and if $\cup_{i=1}^\infty B_i \in {\cal B}_b (\R^d)$, then $L(\cup_{i=1}^\infty B_i) = \sum_{i=1}^\infty L(B_i)$ 
	almost surely,
	\item[(ii)] for every $B \in {\cal B}_b (\R^d)$, the distribution of $L(B)$ is infinitely divisible.
\end{itemize}
We will make a number of assumptions about $g$ and $L$ in the following, which in particular guarantee the existence of the stochastic integral in \eqref{def:X} in the sense of \cite{rajput1989} (see Appendix). 

We assume that for every $B \in {\cal B}_b (\R^d)$, the characteristic function of $L(B)$ has the form
\begin{align}\label{def:L}
\E[ \exp( \i t L(B) ) ] = \exp \Big( \lambda^d (B) \int_{\R_0} ( \exp(\i t y) - 1 - \i t y \1 ( |y| \le 1 ) ) \nu (\d y) \Big), \qquad t \in \R,
\end{align}
where $\nu$ is a symmetric measure on $\R_0 := \R \setminus \{0\}$ satisfying $\int_{\R_0} \min(1,y^2) \nu (\d y) < \infty$.  
Moreover, there exist some $0 \le \beta< 2$, $0<\theta\le 2$ such that
\begin{itemize}
	\item[($\beta$)] 
	$\lim_{y \to 0}  y^\beta\nu ( \{u \in \R_0 : |u|> y\} ) \in (0,\infty)$
	if $\beta >0$, and $\nu (\R_0) < \infty$ if $\beta = 0$,
	\item[($\theta$)] $\limsup_{y \to \infty} y^\theta \nu ( \{u \in \R_0 : |u| > y\} ) < \infty$ if $\theta< 2$, and $\int_{\R_0} y^2 \nu (\d y) < \infty$ if $\theta = 2$,
	\item[($g$)] for every $\boldsymbol{t} \in \R^d$, $g(\boldsymbol{t}, \cdot) \1 (|g(\boldsymbol{t}, \cdot)| \le 1) \in L^\theta (\R^d)$ and   
	$g(\boldsymbol{t},\cdot) \1 (|g(\boldsymbol{t}, \cdot)| > 1) \in L^\beta (\R^d)$.
\end{itemize}
Sometimes we choose $L$ to be a symmetric $\beta$-stable random measure with $0<\beta <2$ and control measure $\lambda^d$, 
i.e.\ for every $B \in {\cal B}_b(\R^d)$, $L(B)$ is a symmetric $\beta$-stable random variable with characteristic function 
$$\E[ \exp(\i t L(B))] = \exp(- \lambda^d (B) |t|^\beta), \qquad t\in \R.$$
In this case the stability index matches the parameter $\beta$ in Assumption~($\beta$) and we can set $\theta=\beta$ in Assumption~($\theta$). 
In the general case the parameter $\beta$ in $(\beta)$ corresponds to  the Blumenthal--Getoor index of $L(B)$:
$$
\beta = \inf \Big\{  q \geq 0 : \int_{0<|y| \le 1} |y|^q \nu (\d y) < \infty \Big\}.
$$
On the other hand, Assumption $(\theta)$ implies that $\int_{|y| > 1} |y|^q \nu (\d y) < \infty$, and hence $\E [|L(B)|^q] < \infty$ for every $0 < q < \theta$ if $\theta< 2$ and $0<q\le \theta$ if $\theta=2$. 

Last, we assume that the kernel $g$ in \eqref{def:X} has the form
\begin{equation}\label{def:g}
	g (\boldsymbol{t}, \boldsymbol{s}) := \sum_{\boldsymbol{\varepsilon} \in \{0,1\}^d} (-1)^{d+\sum_{j=1}^d \varepsilon_j} g_{\boldsymbol{\varepsilon}} (\varepsilon_1 t_1 - s_1, \dots, \varepsilon_d t_d - s_d ), \qquad 
	\boldsymbol{t},\boldsymbol{s} \in \R^d,
\end{equation}
where $g_{\boldsymbol{\varepsilon}} : \R^d \to \R$ is a measurable function for every $\boldsymbol{\vep} \in \{ 0,1 \}^d$. This form of the kernel  is directly motivated by several popular random field models.  Let us present some particular examples. 

\begin{example} \label{ex1} \rm
In cases (ii) and (iv) below $L$ is a symmetric $\beta$-stable random measure with $\beta \in (0,2)$ and control measure $\lambda^d$. \\
(i) A random field $X$ given in \eqref{def:X} is called a L\'evy driven \textit{moving average field} if
\[
g(\boldsymbol{t},\boldsymbol{s}) = g_{(1,\dots,1)} (\boldsymbol{t}-\boldsymbol{s}),
\] 
i.e.\ $g_{\boldsymbol{\varepsilon}} \equiv 0$ for every $\boldsymbol{\varepsilon} \neq (1,\dots,1)$. \\
(ii) It is called a \textit{moving average fractional $\beta$-stable field} (see \cite{take1991}) if 
	$$
	g(\boldsymbol{t},\boldsymbol{s}) = \| \boldsymbol{t}-\boldsymbol{s} \|^{H-\frac{d}{\beta}} - \| \boldsymbol{s} \|^{H-\frac{d}{\beta}}, \qquad  H \in (0,1),~ H \neq \frac{d}{\beta},
	$$
	which corresponds to the choice $g_{(1,\dots,1)} (\boldsymbol{s}) = \| \boldsymbol{s} \|^{H-\frac{d}{\beta}}$, 
	$g_{(0,\dots,0)} (\boldsymbol{s}) = (-1)^{d+1} \| \boldsymbol{s} \|^{H-\frac{d}{\beta}}$ and $g_{\boldsymbol{\varepsilon}}\equiv0$ for every $\boldsymbol{\varepsilon} \neq (1,\dots,1), (0,\dots,0)$. \\
(iii) In \cite{benn2004, cohen2012} a fractional field $X$ has been studied with $\theta=2$ and the kernel
$$
g(\boldsymbol{t},\boldsymbol{s}) = \| \boldsymbol{t}-\boldsymbol{s} \|^{H-\frac{d}{2}} - \| \boldsymbol{s} \|^{H-\frac{d}{2}}, 
\qquad  H \in (0,1),~ H \neq\frac{d}{2}, 
$$ 	
which similarly to the previous example admits the representation \eqref{def:g}. \\
(iv) The \textit{well-balanced symmetric linear fractional $\beta$-stable sheet} $X$ has the kernel
$$
g(\boldsymbol{t},\boldsymbol{s}) = \prod_{i=1}^d ( |t_i - s_i|^{H_i-\frac{1}{\beta}} - |s_i|^{H_i-\frac{1}{\beta}} ), \qquad  H_i \in (0,1),~
H_i \neq \frac{1}{\beta},
$$ 
which can be represented via \eqref{def:g} so that all $g_{\boldsymbol{\varepsilon}}$ are non-trivial. 
Note that $X$ is an extension
of both a well-balanced symmetric linear fractional stable motion, which corresponds to $d=1$, and
of an ordinary fractional Brownian sheet, which corresponds to $\beta=2$. 
\end{example}

\subsection{Power variations and main assumptions} \label{sec2.3}

We consider rectangular increments of the random field $X$ (or any function from $\R^d$ to $\R$) over 
$[\boldsymbol{s},\boldsymbol{t}] = \prod_{i =1}^d [s_i,t_i] \subset \R^d$ for $\boldsymbol{s} < \boldsymbol{t}$, which are defined as
\begin{align}\label{def:incrX}
X ( [\boldsymbol{s},\boldsymbol{t}] )
&:= \sum_{\boldsymbol{\varepsilon} \in \{0,1\}^d} (-1)^{d+\sum_{j=1}^{d} \varepsilon_j} X (s_1+\varepsilon_1 (t_1-s_1), \dots,s_d+\varepsilon_d (t_d-s_d)). 
\end{align}
For instance, when $d=1$ \eqref{def:incrX} reduces to $X([s,t]) = X(t) - X(s)$, while 
$
X([\boldsymbol{s},\boldsymbol{t}]) = X (t_1,t_2) - X (t_1,s_2) - X (s_1,t_2) + X (s_1,s_2)
$  when  $d = 2$. 
The rectangular increment can also be recovered by differencing iteratively with respect to each of the arguments of $X$, that is 
$$
X([\boldsymbol{s},\boldsymbol{t}]) =  \Delta_{t_1-s_1}^{(1)} \dots \Delta_{t_d-s_d}^{(d)} X (\boldsymbol{s}),
$$
where
$
\Delta_{t_i-s_i}^{(i)} X (\boldsymbol{s}) = X(\boldsymbol{s} + (t_i-s_i) \boldsymbol{e}_i) - X (\boldsymbol{s})
$
is a directional increment, $i = 1, \dots, d$, and $\{ \boldsymbol{e}_1, \dots, \boldsymbol{e}_d \}$ is the standard basis of $\R^d$. 
The random field $X$ in \eqref{def:X} has stationary rectangular increments, i.e.\ for any fixed $\boldsymbol{s} \in \R^d$, 
$$
(X([\boldsymbol{s},\boldsymbol{t}]))_{ \boldsymbol{s} < \boldsymbol{t}} \overset{\textnormal{fdd}}{=} ( X([\boldsymbol{0},\boldsymbol{t}-\boldsymbol{s}]))_{\boldsymbol{s} < \boldsymbol{t}}.
$$ 
Indeed, the rectangular increment  of the function $g(\cdot,\boldsymbol{u})$ in \eqref{def:g} over $[\boldsymbol{s},\boldsymbol{t}]$ coincides with that of $g_{(1,\dots,1)}$ over $[\boldsymbol{s}-\boldsymbol{u},\boldsymbol{t}-\boldsymbol{u}]$, while  all of the other functions $g_{\boldsymbol{\vep}}$, 
$\boldsymbol{\vep} \neq (1,\ldots, 1)$, vanish after the computation of the rectangular increments (but 
they are usually still needed for 
the stochastic integrals in \eqref{def:X} to exist).   Since only the function $g_{(1,\dots,1)}$ matters when taking rectangular increments, we write with a slight abuse of notation  
\begin{equation}\label{ma}
g(\boldsymbol{s}) = g_{(1,\dots,1)}(\boldsymbol{s}), \qquad  \boldsymbol{s} \in \R^d.
\end{equation}
We also write $\Delta_r X(\boldsymbol{s})$ for $X ([\boldsymbol{s}, \boldsymbol{s}+r \boldsymbol{1}])$, where $\boldsymbol{1}=(1,\ldots, 1)\in \R^d$ and all edges of the rectangle have equal length $r>0$.

Our main focus are power variation statistics of $X$ computed over the set $[0,1]^d$: 
\begin{equation}\label{def:Vpn}
V_n(p) := \sum_{\boldsymbol{i} \in \{0,\dots,n-1\}^d
} | \Delta_{1/n} X(\boldsymbol{i}/n) |^p
\end{equation} 
for $p>0$. The main goal of this paper is to study the asymptotic behaviour of the statistic $V_n(p)$ as $n \to \infty$. We will see that the type and mode of the limit crucially depend on the behaviour of the function $g: \R^d \to \R$ introduced in \eqref{ma}. More specifically, we will assume that $g$ is \textit{locally homogenous} near $\boldsymbol{0}$. That is, we consider
$g (\boldsymbol{s}) \sim h(\boldsymbol{s})$
as $\boldsymbol{s} \to \boldsymbol{0}$,
where $h$ is an absolutely homogeneous function of some degree $\delta \neq 0$, i.e.\ $h (a \boldsymbol{s}) = |a|^{\delta} h (\boldsymbol{s})$ for all $a \in \R$ and $\boldsymbol{s} \in \R^d$. However, this type of assumption still does not uniquely determine the asymptotic theory in contrast to the 
theory of case $d=1$ investigated in \cite{basse2017}. We will therefore distinguish two classes of homogeneous functions $h: \R^d \to \R$:

\begin{itemize}
	\item[(H1)] For all $\boldsymbol{s} \in \R^d$, 
	\begin{equation*}
	g (\boldsymbol{s}) = f (\boldsymbol{s}) h (\boldsymbol{s}), \qquad \text{where }
	h (\boldsymbol{s}) := \| \boldsymbol{s} \|^{d\alpha} \text{ for some } 
	\alpha \neq 0,
	\end{equation*}
	and $f$  has continuous partial derivatives up to the $d$-th order at every point in $\R^d$ and  $f(\boldsymbol{0}) = 1$. Moreover, there exists $\rho > 0$ such that $|\partial^d g|$ is in $L^\theta ( B_\rho^c (\boldsymbol{0}) )$ and is radially non-increasing, i.e.\ $|\partial^d g (\boldsymbol{s})| \ge |\partial^d g (\boldsymbol{t})|$ if $\rho \le \|\boldsymbol{s}\| \le \|\boldsymbol{t}\|$, $\boldsymbol{s}, \boldsymbol{t} \in \R^d$.
	
	\item[(H2)] For all $\boldsymbol{s} \in \R^d$, $g (\boldsymbol{s}) = \prod_{i =1}^d g_i (s_i)$.
	For all $s \in \R$, 
	$$
	g_i (s)=  f_i (s) h_i (s), \qquad \text{where } h_i (s) := |s|^{\alpha_i} 
	\text{ for some } 
	\alpha_i \neq 0,
	$$
	and $f_i \in C^1 (\R)$ satisfies $f_i (0) = 1$, $i=1,\dots,d$. Moreover, there exists $\rho > 0$ such that $g'_i \in L^q ( (-\rho,\rho)^c )$ with $q := \min(\theta,\max(\beta,p))$ and
	$|g'_i(s)| \ge |g'_i(t)|$ if $\rho \le |s| \le |t|$, $s,t \in \R$, $i=1,\dots,d$. We set 
	$$f(\boldsymbol{s}) := \prod_{i =1}^d f_i (s_i), \qquad h(\boldsymbol{s}) := \prod_{i=1}^d h_i (s_i), \qquad \boldsymbol{s} \in \R^d.$$ 
\end{itemize}

\noindent
We will see in the next section that under (H1), where the homogeneous function $h$ does not depend on the direction, the limit theory for the power variation $V_n(p)$ in some sense resembles the case $d=1$ studied in
\cite{basse2017}.  On the other hand, the asymptotic results for kernel satisfying (H2) 
	are more complex because they allow for mixtures in terms of conditions and limits obtained before.

\begin{remark} \rm
The assumption $f(\boldsymbol{0}) = 1$  in (H1) is not essential (the same applies to the corresponding assumption in (H2)). 
As long as $f(\boldsymbol{0}) \neq 0$ we may deduce the setting of (H1) by adjusting the L\'evy measure $\nu$ accordingly. In (H2) the multiplicative form of the homogeneous function $h$ is essential, while the analogous assumption on the function $f$ is not necessary and it
is considered for simplicity of exposition. 
\qed
\end{remark}

\section{Main results} \label{sec3}
\setcounter{equation}{0}
\renewcommand{\theequation}{\thesection.\arabic{equation}}

In this section we consider the random field  $X=(X(\boldsymbol{t}))_{\boldsymbol{t}\in\R^d}$ defined  in \eqref{def:X} with $L$ and $g$ given by \eqref{def:L} and \eqref{def:g}, respectively, and satisfying  Assumptions ($g$), ($\theta$) and  ($\beta$)  for some $0<\theta \le 2$ and  $0 \le \beta < 2$. The  two following theorems state   the limit theory for power variation statistics $V_n(p)$ of $X$ under (H1) and (H2).   
Its mode of convergence and limit depend on the interplay between the power $p$, the Blumenthal--Getoor index $\beta$
and the form of the kernel $g$ at the origin. In each case we use the most convenient representation of $X$ or $L$. In Theorem \ref{thm1}(i) we  will use 
a Poisson random measure $\Lambda^\dagger$ on $[0,1]^d \times \R_0$ with intensity measure $\lambda^d \otimes \nu$, which is constructed by adding to 
the  jump sizes of $L$ restricted to  $[0,1]^d$, the marks that are 
i.i.d.\ random vectors with a common uniform distribution on $[0,1]^d$, defined on the extension of the underlying probability space $(\Omega, {\cal F}, \P)$ and independent of the $\sigma$-algebra ${\cal F}$. Similarly, in Theorem \ref{thm2}(i) a Poisson random measure $\Lambda^\ddagger$ with intensity measure $\lambda^k \otimes \lambda^{d-k} \otimes \nu$ on $[0,1]^k \times \R^{d-k} \times \R_0$ is constructed from the  jumps of $L$ on $[0,1]^k \times \R^{d-k}$ for some $k=1,\dots,d$. First we state the limit  theory for the statistic $V_n(p)$ under (H1).

\begin{theorem}\label{thm1}
Let Assumption \textnormal{(H1)} hold for some $\alpha\in \R_0$. 
	
	\begin{itemize}
		\item[(i)] 
		Let $p > \beta$ and $\alpha + 1/p \in (0,1)$. Then 
		\begin{equation*}
		n^{d \alpha p} V_n(p) \overset{{\cal F}\textnormal{-d}}{\to} \int_{[0,1]^d \times \R_0} \Big( | y |^p \sum_{\boldsymbol{j} \in \Z^d} |\Delta_1 h ( \boldsymbol{j} - \boldsymbol{u}) |^p \Big) \Lambda^\dagger (\d \boldsymbol{u}, \d y)
		\end{equation*}
		where $\Lambda^\dagger$ is the  Poisson random measure on $[0,1]^d \times \R_0$ having intensity measure $\lambda^d \otimes \nu$ defined in Definition~\ref{ljsdlfjsd}.
		
		\item[(ii)] Let $L$ be a symmetric $\beta$-stable random measure on $\R^d$ with $\beta \in (0,2)$ and control measure $\lambda^d$. 
		Let $p < \beta = \theta$ and 
		$H := \alpha +1/\beta \in (0,1)$. 
		Then
		\begin{equation*}
		n^{d(Hp-1)} V_n(p) \overset{L^1}{\to} \E[ | L( [0,1]^d ) |^p ]\Big( \int_{\R^d} | \Delta_1 h( \boldsymbol{s}) |^\beta \d \boldsymbol{s} \Big)^{\frac{p}{\beta}}.
		\end{equation*}
		
		\item[(iii)] Let $p \ge 1$ and $\alpha + 1/\max(\beta, p) > 1$. Then 
		$$
		n^{d(p-1)} V_n(p) \overset{\textnormal{a.s.}}{\to}
		\int_{[0,1]^d} | Y (\boldsymbol{t}) |^p \d \boldsymbol{t}, 
		$$		
		where $(Y(\boldsymbol{t}) )_{\boldsymbol{t} \in [0,1]^d}$ is a measurable random field satisfying
		$$
		Y (\boldsymbol{t}) = \int_{\R^d} \partial^d g(\boldsymbol{t} - \boldsymbol{s}) L(\d \boldsymbol{s})\text{ a.s.} \qquad \text{for all }\boldsymbol{t}\in [0,1]^d, 
		 	$$ 
		and 
		$$
		\int_{[0,1]^d} |Y(\boldsymbol{t})|^p \d \boldsymbol{t} < \infty \text{ a.s.}
		$$
			\end{itemize}
\end{theorem}

\noindent

We note  that Theorem~\ref{thm1} covers all $\alpha\in \R_0$ satisfying  $\alpha>-1/\max(\beta,p)$ except for the three boundary cases 
$p=\beta$ and  $\alpha=1-1/\max(\beta,p)$  with the additional assumption that $L$ is $\beta$-stable 
if both $p<\beta$ and $\alpha<1-1/\beta$, and  with the additional assumption that $p\geq 1$ if $\alpha>1-1/\max(\beta,p)$. Remark that we obtain very different convergence rates and types/modes of limits in Theorem~\ref{thm1}. 
While Theorem \ref{thm1}(ii) is 
of ergodic type, Theorem \ref{thm1}(i) and (iii) are quite non-standard. A similar phenomenon has been observed 
for processes in \cite{basse2017}. Indeed, the results of Theorem \ref{thm1} look like a direct extension of \cite[Theorem 1.1]{basse2017} from $d=1$
to a general dimension $d\geq 1$. In contrast to the imposed assumptions in \cite[Theorem 1.1]{basse2017}, (H1) allows for negative values of $\alpha$.
The next result presents the asymptotic theory for the statistic $V_n(p)$ under (H2).

\begin{theorem}\label{thm2}
	Let Assumption \textnormal{(H2)} hold for some $\alpha_1,\dots,\alpha_d\in \R_0$, and $p\neq \theta$ if $\theta <2$. 
		
	\begin{itemize}
\item[(i)] Let $p > \beta$. 
For some $k=1,\dots,d$ let $\alpha_i +1/p \in (0,1)$ for  $i=1,\dots,k$, and $\alpha_i+1/p>1$ for $i=k+1,\dots,d$. Then 
\begin{align*}
& n^{(d-k) (p-1) +\sum_{i=1}^{k} \alpha_i p} V_n(p) \overset{{\cal F}{\textnormal{-d}}}{\to}
\\ &\qquad 
\int_{[0,1]^k \times \R^{d-k} \times \R_0} \Big( |y|^p \Big( \prod_{i=1}^k \sum_{j \in \Z} |\Delta_1 h_i(j-u_i)|^p\Big) \prod_{i =k+1}^{d} \int_0^1 |g'_i (t-x_i) |^p \d t \Big) \Lambda^\ddagger (\d \boldsymbol{u},\d \boldsymbol{x}, \d y),
\end{align*}
where $\Lambda^\ddagger$ is the Poisson random measure on $[0,1]^k \times \R^{d-k} \times \R_0$ having intensity measure $\lambda^k \otimes \lambda^{d-k} \otimes \nu$ defined in Definition~\ref{sldjfsghsklhs}, and $\boldsymbol{u}=(u_1,\dots,u_k) \in [0,1]^k$, $\boldsymbol{x}=(x_{k+1},\dots,x_{d}) \in \R^{d-k}$. 

\item[(ii)]  Let $L$ be a symmetric $\beta$-stable random measure on $\R^d$ with $\beta \in (0,2)$ and control measure $\lambda^d$. 
Let $p < \beta = \theta$. For some $k =1,\dots,d$ let
$H_i:=\alpha_i + 1/\beta \in (0,1)$ for  $i=1,\dots,k$, and $\alpha_i + 1/\beta > 1$ for  $i=k+1,\dots,d$. Then
\begin{equation}\label{lim2}
n^{(d-k)(p-1)+\sum_{i=1}^k (H_ip-1)} V_n(p) \overset{L^1}{\to} \E[ | L([0,1]^d) |^p] \prod_{i =1}^k \Big( \int_{\R} |\Delta_1 h_i (s)|^\beta \d s \Big)^{\frac{p}{\beta}} \prod_{i=k+1}^d \Big( \int_{\R} |g_i'(s)|^\beta \d s \Big)^{\frac{p}{\beta}}.
\end{equation}

\item[(iii)] Let $p \ge 1$ and $\alpha_i+ 1/\max(\beta, p) > 1$, $i=1,\dots,d$. 
Then
$$
n^{d(p-1)} V_n(p) \overset{\textnormal{a.s.}}{\to} \int_{[0,1]^d} | Y (\boldsymbol{t}) |^p \d \boldsymbol{t},
$$
where $( Y(\boldsymbol{t}) )_{\boldsymbol{t} \in [0,1]^d}$ is a measurable random field satisfying 
\begin{equation}\label{ljsdlfjhsgdhls}
 Y (\boldsymbol{t})= \int_{\R^d} \prod_{i =1}^d g'_i (t_i - s_i) L(\d \boldsymbol{s}) \text{ a.s.} \qquad \text{for all } \boldsymbol{t} \in [0,1]^d,
\end{equation}
and
$$
\int_{[0,1]^d} |Y(\boldsymbol{t})|^p \d \boldsymbol{t} < \infty \text{ a.s. }
$$
\end{itemize}
\end{theorem}

\noindent
Under Assumption \textnormal{(H2)} there is no loss of generality by assuming that $\alpha_1\leq \alpha_2\leq \dots\leq \alpha_d$, and therefore Theorem~\ref{thm2} covers all $\alpha_1,\dots,\alpha_d\in \R_0$ with $\alpha_1>-1/\max(\beta,p)$ except for the boundary cases where $p=\beta$  or $\alpha_k=1-1/\max(\beta,p)$ for some $k=1,\dots,d$ with the two additional assumptions that $L$ is $\beta$-stable if both $p<\beta$  and $\alpha_k+1/\beta<1$ for some $k=1,\dots,d$, and moreover that  $p\geq 1$ if  $\alpha_i+1/\max(\beta, p)>1$ for all $i=1,\dots,d$.  
The results of Theorem~\ref{thm2} are more complex compared to the isotropic type setting of Theorem \ref{thm1}. Since we have more degrees
of freedom for the powers $\alpha_i$ under Assumption (H2) than under Assumption (H1),  certain mixtures of  Theorem \ref{thm1}(i)--(iii) appear 
in Theorem \ref{thm2}. Indeed, when $p>\beta$ and the first $k$ indices $\alpha_i$ satisfy the assumption of  Theorem \ref{thm1}(i) while the last 
ones satisfy the assumption of  Theorem \ref{thm1}(iii), we obtain their mixture in Theorem \ref{thm2}(i). Similarly, Theorem \ref{thm2}(ii) can be interpreted as a mixture of    Theorem \ref{thm1}(ii) and (iii).

\begin{remark} \rm \label{remclt}
(i) Theorems \ref{thm1}(ii) and \ref{thm2}(ii) remain valid for $\beta =2$, where $L$ is a Gaussian random measure on $\R^d$ with  zero mean and variance $\lambda^d$. In this case the result holds true for all $p > 0$. \\
(ii) Assume that the function $h$ satisfies (H2) with $\alpha_1 = \dots = \alpha_d$. Then 
we have $k = d$  in Theorem~\ref{thm2}(i) and (ii). 
Furthermore, rates of convergence and limits of $V_n(p)$ coincide with those in Theorem \ref{thm1}, which implies that we cannot distinguish 
between the classes (H1) and (H2) based upon the statistic $V_n(p) $. \qed
\end{remark}

\noindent
Next, we examine how the results of Theorems \ref{thm1} and \ref{thm2} apply to models discussed in Example \ref{ex1}.

\begin{example} \label{ex2} \rm (Continuation of Example \ref{ex1})
In cases (ii), (iv) and (v) below let $L$ be a symmetric $\beta$-stable random measure with $\beta \in (0,2)$ and control measure $\lambda^d$. In all cases let $p>0$.\\
(i): We consider a special case of a L\'evy driven moving average field $X$ having
$$
g(\boldsymbol{t},\boldsymbol{s}) = g_{(1,\dots,1)}(\boldsymbol{t}-\boldsymbol{s}) \qquad \text{with } g_{(1,\dots,1)}(\boldsymbol{s}) = \frac{2}{\Gamma (\frac{d}{4}-\frac{\gamma}{2})} \Big\| \frac{2\boldsymbol{s}}{\sigma}  \Big\|^{\frac{\gamma}{2}-\frac{d}{4}} K_{\frac{\gamma}{2}-\frac{d}{4}} ( \sigma \| \boldsymbol{s}\| ),
$$
where $\gamma \in (0, d/2)$, $\sigma>0$ and $K_{\gamma/2-d/4}$ denotes the modified Bessel function of the second kind.
It holds that 
$$K_{\frac{\gamma}{2}-\frac{d}{4}}(s) \sim \frac{1}{2} \Gamma\Big(\frac{d}{4}-\frac{\gamma}{2}\Big) \Big(\frac{s}{2}\Big)^{-(\frac{d}{4}-\frac{\gamma}{2})}
\qquad \text{as }s \downarrow 0,
$$ 
see \cite[Eq.\ (9.6.9), p.\ 375]{abra1972}. This  implies
$g_{(1,\dots,1)}(\boldsymbol{s}) \sim  \| \boldsymbol{s} \|^{\gamma-\frac{d}{2}}$ as $\boldsymbol{s} \to \boldsymbol{0}$. It has been shown in \cite{hansen2013,jons2013} that such a choice of $g$ induces a covariance function
\begin{align*}
\operatorname{Cov}(X(\boldsymbol{0}), X(\boldsymbol{t})) = \operatorname{Var} (X (\boldsymbol{0})) \frac{2^{1-\gamma}}{\Gamma (\gamma)} ( \sigma \| \boldsymbol{t} \| )^\gamma K_\gamma (\sigma \| \boldsymbol{t} \|), \qquad \boldsymbol{t} \in \R^d,
\end{align*}
belonging to the \textit{Mat\'ern family} when $\E[X (\boldsymbol{0})^2]<\infty$  (see \cite{gutt2006} for more details). Then Theorem~\ref{thm1}(i) applies if $p> \beta$ and $(1/2-1/p)d<\gamma<(3/2-1/p)d$,
 Theorem~\ref{thm1}(ii) applies if $p< \beta$, $(1/2-1/\beta)d<\gamma<(3/2-1/\beta)d$ and $L$ is $\beta$-stable. Theorem~\ref{thm1}(iii)
 never applies to this example.
 \\
(ii): The kernel 
$$
	g(\boldsymbol{t},\boldsymbol{s}) = \| \boldsymbol{t}-\boldsymbol{s} \|^{H-\frac{d}{\beta}} - \| \boldsymbol{s} \|^{H-\frac{d}{\beta}}, \qquad  H \in (0,1),~ H \neq \frac{d}{\beta},
$$
satisfies (H1). Therefore,  Theorem~\ref{thm1}(i) applies if $p>\beta$ and $H>(1/\beta-1/p)d$,  
  Theorem~\ref{thm1}(ii) applies if   $p<\beta$. Again Theorem~\ref{thm1}(iii) never applies  for this example.  \\
(iii): The kernel
$$
g(\boldsymbol{t},\boldsymbol{s}) = \| \boldsymbol{t}-\boldsymbol{s} \|^{H-\frac{d}{2}} - \| \boldsymbol{s} \|^{H-\frac{d}{2}}, 
\qquad  H \in (0,1),~ H \neq\frac{d}{2},
$$ 	 
obviously satisfying (H1), induces the covariance function 
$$
\operatorname{Cov} (X(\boldsymbol{t}), X(\boldsymbol{s})) = \operatorname{Var}( X(\boldsymbol{e}_1) )  \frac{1}{2} ( \| \boldsymbol{s} \|^{2H} + \|\boldsymbol{t}\|^{2H}-\| \boldsymbol{t}-\boldsymbol{s}\|^{2H}  ), \qquad \boldsymbol{t},\boldsymbol{s} \in \R^d,
$$
when $\E[X(\boldsymbol{e}_1)^2]<\infty$.
Hence Theorem~\ref{thm1}(i) applies if $p>\beta$ and $(1/2-1/p)d<H<(3/2-1/p)d$,  Theorem~\ref{thm1}(ii) applies  if $p<\beta$, $H<(3/2-1/\beta)d$ and $L$ is $\beta$-stable. Theorem~\ref{thm1}(iii) never applies to this example.
\\
(iv): The kernel
$$
g(\boldsymbol{t},\boldsymbol{s}) = \prod_{i=1}^d ( |t_i - s_i|^{H_i-\frac{1}{\beta}} - |s_i|^{H_i-\frac{1}{\beta}} ), \qquad  H_i \in (0,1),~
H_i \neq \frac{1}{\beta},
$$ 
satisfies assumption (H2)  with $\alpha_i =H_i-1/\beta$, $i=1,\dots,d$, and $q=\beta$. We may and do assume that $H_1\leq H_2\leq \dots\leq H_d$. Therefore, Theorem~\ref{thm2}(i) applies if $H_1>1/\beta-1/p$
and $p>\beta$,  Theorem~\ref{thm2}(ii) applies if  $p< \beta$, whereas Theorem~\ref{thm2}(iii) never applies to this example. \\
(v): Recalling the notation of rectangular increments we introduce a new kernel 
$$
g(\boldsymbol{t},\boldsymbol{s}) = h ( [-\boldsymbol{s},\boldsymbol{t}-\boldsymbol{s}]) \qquad \text{with } h (\boldsymbol{s}) = \| \boldsymbol{s} \|^{d(H-\frac{1}{\beta})}, \qquad H \in (0,1),~ H \neq \frac{1}{\beta}.
$$
In particular, when $d=2$ it holds that 
$g(\boldsymbol{t},\boldsymbol{s}) = h(t_1-s_1,t_2-s_2)-h(t_1-s_1,-s_2)-h(-s_1,t_2-s_2)+h(-s_1,-s_2)$.
In this case (H1) is satisfied and 
Theorem~\ref{thm1}(i) applies if $H>1/\beta-1/p$ and $p>\beta$,  Theorem~\ref{thm1}(ii) applies if $p<\beta$. Theorem~\ref{thm1}(iii) never applies to this example.  \qed
\end{example}

\noindent
Theorems \ref{thm1} and \ref{thm2} have important consequences for parameter identification and parameter estimation. To illustrate the potential 
of Theorem \ref{thm1} let us consider the moving average fractional $\beta$-stable field defined in Example \ref{ex1}(ii). A standard strategy to estimate the Hurst parameter $H \in (0,1)$ is to use a ratio statistic based on a change of frequency. More specifically, the ergodic result of Theorem \ref{thm1}(ii) immediately implies the convergence
\[
R_n := \frac{ \sum_{\boldsymbol{i}\in \{0,\dots,n-2\}^d
	} | \Delta_{2/n} X(\boldsymbol{i}/n) |^p}{ \sum_{\boldsymbol{i} \in \{0,\dots,n-1\}^d
	} | \Delta_{1/n} X(\boldsymbol{i}/n) |^p} \overset{\P}{\to} 2^{dHp}
\]
if $p<\beta$. Hence, 
\begin{equation}\label{def:estimator}
H_n:= \frac{\log R_n}{dp\log 2}\overset{\P}{\to} H, \qquad \text{if } p<\beta.
\end{equation}
Obviously, the proposed estimation procedure assumes prior knowledge of the parameter $\beta$, since we need to choose $p\in (0,\beta)$.  In 
the case $d=1$ the papers 
\cite{DI17,LP19,LP20,MOP20} have suggested to use negative powers $p\in (-1,0)$ to estimate the parameter $H$ for unknown $\beta$. A similar idea should apply in the random field setting, although negative power variations are beyond the scope of our paper. A construction of confidence regions for parameters of the moving average fractional $\beta$-stable field requires proving the weak limit theory associated with Theorem 
\ref{thm1}(ii). However, this is a rather complex problem since the martingale type techniques, which have been applied 
for processes in \cite{BHP19,basse2017}, do not easily extend to our setting.  

A straightforward consequence of Theorems \ref{thm1} and \ref{thm2} is the identification of some involved parameters via the corresponding convergence rates. Indeed, we observe that the statistic $S_n(p):=\log V_n(p)/ \log n$ converges in probability to the exponent of the convergence rates given in Theorems \ref{thm1} and \ref{thm2}. Considering again the moving average fractional $\beta$-stable field as an example, the three convergence rates described in  
Theorem~\ref{thm1} and the points of phase transition uniquely determine the parameter $(H,\beta)$. In other words,
the limit of the process $(S_n(p))_{p> 0}$ identifies  $(H,\beta)$. The same logic applies to the well-balanced symmetric linear fractional $\beta$-stable sheet discussed in Example \ref{ex1}(iv), where the limit of  $(S_n(p))_{p> 0}$ uniquely determines the 
parameter $(\sum_{i=1}^d H_i, \beta)$; however, Theorem \ref{thm2} does not suffice to identify/estimate the parameters $(H_i)_{1\leq i\leq d}$ separately. To provide such an inference we can identify/estimate
 $H:=H_i$ from increments of a line process $(X(\boldsymbol{1} + t \boldsymbol{e}_i))_{t \in \R}$. Indeed, it is a well-balanced symmetric linear $\beta$-stable motion, to which  Theorem \ref{thm1}(ii) applies with $d=1$. Hence, we may obtain a consistent estimator of $H_i$ via \eqref{def:estimator}.

%



\section{Proofs} \label{sec4}
\setcounter{equation}{0}
\renewcommand{\theequation}{\thesection.\arabic{equation}}



We first present some preliminary facts that will be used in the proofs. We will use a stable convergence of fractional parts of random variables: if $\boldsymbol{W} \sim {\cal U}([0,1]^d)$, then as $n \to \infty$,
\begin{equation} \label{braconv}
\{ n \boldsymbol{W} \} \overset{{\cal F}{\textnormal{-d}}}{\to} \boldsymbol{U},
\end{equation}
where $\boldsymbol{U}$ is ${\cal U}([0,1]^d)$-distributed random vector, defined on the extension 
of the underlying probability space $(\Omega,{\cal F},\P)$ and independent of the $\sigma$-algebra ${\cal F}$; see e.g. \cite[Lemma 4.1]{basse2017}. We will repeatedly use the following inequalities. Let $m \in \N$, $p>0$. For $\boldsymbol{a} \in \R^m$,  set $\| \boldsymbol{a} \|_p = (\sum_{i=1}^m |a_i|^p )^{1/p}$. For $\boldsymbol{a}, \boldsymbol{b} \in \R^m$, it holds that
\begin{align}
| \| \boldsymbol{a} \|_p^p - \| \boldsymbol{b} \|_p^p  | &\le \| \boldsymbol{a} - \boldsymbol{b} \|_p^p \qquad \text{if } 0 < p \le 1,
\label{ineq0}\\
| \| \boldsymbol{a} \|_p - \| \boldsymbol{b} \|_p  | &\le \| \boldsymbol{a} - \boldsymbol{b} \|_p \qquad \text{if } p > 1.
\label{ineq1}
\end{align}
For an  $n\in \N$ we  set $\boldsymbol{n} := (n,\dots,n) \in \N^d$.

\subsection{Some Poisson random measures related to $L$}\label{sec:PRM}

By extending our probability space $(\Omega, \mathcal F,\P)$ if necessary we  may and do assume that it is rich enough to support a $\mathcal U([0,1])$-distributed random variable independent of $L$. 
To the infinitely divisible random measure $L$ given in \eqref{def:L}, we associate a random field $(L(\boldsymbol{t}))_{\boldsymbol{t}\in\R^d}$ by $L(\boldsymbol{t}) = L([\boldsymbol{0},\boldsymbol{t}])$ (for $\boldsymbol{t} > \boldsymbol{0}$, and similarly otherwise). We note that $(L(\boldsymbol{t}))_{\boldsymbol{t}\in \R^d}$ is a L\'evy process in the sense of \cite[page~5]{AMSW} and in particular for all $n\in \N$ and all disjoint rectangles
$[\boldsymbol{a}_1,\boldsymbol{b}_1], \dots, [\boldsymbol{a}_n,\boldsymbol{b}_n]$ in $\R^d$,  $L([\boldsymbol{a}_1,\boldsymbol{b}_1]),\dots, L([\boldsymbol{a}_n,\boldsymbol{b}_n])$ are independent. As c\`adl\`ag functions of several variables are less standard than the univariable case, we will define the appropriate sample path space for $(L(\boldsymbol{t}))_{\boldsymbol{t} \in \R^d}$ in the following. For $d=1,2,\dots$ we say that a function $x\!:\R^d \to \R$ is \emph{lamp} (limits along monotone paths) if for all $\boldsymbol{t}\in \R^d$ we have 
\begin{enumerate} 
	\item    the limit $x(\boldsymbol{t},\mathcal R):=\lim_{\boldsymbol{u}\to \boldsymbol{t},\, \boldsymbol{u}\mathcal R \boldsymbol{t}} x(\boldsymbol{u})$ exists in $\R$ for each of the $2^d$ order relations $\mathcal R=(R_1,\dots,R_d)$, where $R_i$ is either $\geq $ or $<$ for $i=1,\dots, d$,
	\item $x(\boldsymbol{t})=x(\boldsymbol{t},\mathcal R)$ when $\mathcal R=(\geq,\dots, \geq )$.  
\end{enumerate}
For each lamp function $x\!\!:\R^d \to \R$ we define the point mass  jump $J_{\boldsymbol{t}}(x)$ of $x$ at $\boldsymbol{t}\in \R^d$ as 
$J_{\boldsymbol{t}}(x)=\lim_{\boldsymbol{u} \to \boldsymbol{t},\, \boldsymbol{u}\mathcal R \boldsymbol{t}} x([\boldsymbol{u},\boldsymbol{t}])$, where $\mathcal R=(<,\dots,<)$.  For instance, when $d=1$, we have  $J_t(x)=x(t)-x(t-)$, where  $x(t-)=x(t,<)$ denotes the left-hand limit, while $J_{\boldsymbol{t}}(x)=x(t_1,t_2)-x(t_1,t_2-)-x(t_1-,t_2)+x(t_1-,t_2-)$ when $d=2$.  
The above notation and terminology are due to Straf~\cite{Straf}. By Proposition~4.1 of \cite{AMSW} and homogeneity of $L$, $( L(\boldsymbol{t}))_{\boldsymbol{t}\in \R^d}$ has a lamp modification, which also will be denoted $( L(\boldsymbol{t}))_{\boldsymbol{t}\in \R^d}$. For every Borel set $A$ of $\R^d\times \R_0$, set 
\begin{equation}\label{ljsdlfjsdlj}
\Lambda(A)= \# \{\boldsymbol{v}\in \R^d: (\boldsymbol{v},J_{\boldsymbol{v}}(L))\in A\}, 
\end{equation}
where $\# S$ denotes the number of elements in a set $S$. 
From Proposition~4.4 of \cite{AMSW} we deduce that $\Lambda$ is a Poisson random measure on $\R^d \times \R_0$ with intensity measure $\lambda^d\otimes \nu$, and by Theorem~4.6 of \cite{AMSW} we have that for all $\boldsymbol{t} \in \R^d$,
\begin{align*}
L(\boldsymbol{t})= {}&  \int_{(\boldsymbol{0},\boldsymbol{t}]\times \{|y|>1\}} y\,\Lambda( \d\boldsymbol{v},\d y)+ \lim_{\epsilon\downarrow 0} \int_{(\boldsymbol{0},\boldsymbol{t}]\times \{\epsilon<|y|\leq 1\}} y \, \big( \Lambda(\d\boldsymbol{v},\d y)- (\lambda^d\otimes \nu) (\d \boldsymbol{v},\d y) \big) \\
={}& \lim_{\epsilon\downarrow 0}  \int_{(\boldsymbol{0},\boldsymbol{t}]\times \{|y|>\epsilon\}} y\,\Lambda(\d \boldsymbol{v},\d y) =: \int_{(\boldsymbol{0},\boldsymbol{t}]\times \R_0} y\,\Lambda(\d \boldsymbol{v},\d y)
\end{align*}
where the second equality follows by symmetry of $\nu$, and the convergence to the two limits is uniform in $\boldsymbol{t}$ on compact subsets of $\R^d$ almost surely.

In the following we will construct a proper point process representation of $\Lambda$ restricted to $[0,1]^d\times \R_0$,
which we are going to use in Theorem~\ref{thm1}(i).  Since $\nu$ is a $\sigma$-finite measure we may choose a probability measure $\tilde \nu$ such that $\nu$ is absolute continuous with respect to $\tilde \nu$ with density $\rho>0$. Let  $(W_k)_{k\in \N}$ be  an i.i.d.\ sequence of real-valued random variables with the common distribution $\tilde \nu$, 
$(\tilde {\boldsymbol{V}}_k)_{k\in \N}$ be an i.i.d.\ sequence of  $\mathcal U([0,1]^d)$-distributed random vectors, and 
$(\Gamma_k)_{k\in \N}$ be  a sequence of partial sums of i.i.d.\ standard exponential  random variables. 
Assume that the three sequences $(\tilde {\boldsymbol{V}}_k)_{k\in \N}$, $(W_k)_{k\in \N}$ and $(\Gamma_k)_{k\in \N}$ 
are independent, and set 
\begin{equation}
\tilde J_k= W_k \1(\rho(W_k)\geq \Gamma_k),\  k\in \N, \qquad \text{and }\qquad \tilde \Lambda=\sum_{k=1}^\infty \delta_{(\tilde {\boldsymbol{V}}_k,\tilde J_k)}.
\end{equation}
Then $\tilde \Lambda$ is a Poisson random measure on $[0,1]^d \times \R_0$
with intensity measure $\lambda^d\otimes \nu$, and since our probability space supports an $\mathcal U([0,1])$-random variable independent of $L$ by assumption, there exists a sequence $({\boldsymbol{V}}_k,J_k)_{k\in \N}$ which equals $(\tilde {\boldsymbol{V}}_k,\tilde J_k)_{k\in \N}$ in distribution, and satisfies 
\begin{equation}\label{sdljflsjd}
\Lambda 
=\sum_{k=1}^\infty \delta_{({\boldsymbol{V}}_k,J_k)}
\end{equation}
on $[0,1]^d \times \R_0$ almost surely,
cf.\  Proposition~2.1 in \cite{Ros-Point}. 
In the following we will describe some Poisson random measures 
appearing in the limit of  Theorem~\ref{thm1}(i).

\begin{definition}\label{ljsdlfjsd}
	Let $(\boldsymbol{U}_k)_{k\in \N}$ be an i.i.d.\ sequence of $\mathcal U([0,1]^d)$-distributed random vectors, defined on an extension of $(\Omega,\mathcal F,\P)$ and independent of $\mathcal F$, and set 
	\begin{equation}\label{ljsdfljshslp}
	\Lambda^\dagger =\sum_{k=1}^\infty \delta_{(\boldsymbol{U}_k,J_k)}. 
	\end{equation}
\end{definition}

For Theorem~\ref{thm2}(i) we need a proper point process representation of $\Lambda$ restricted 
to $[0,1]^k\times \R^{d-k}\times \R_0$, where $k=1,\dots,d$. 
To this aim, let us introduce a probability measure $\kappa$ on $\R^{d-k}$ by $\kappa(\d \boldsymbol{x})=h_1(\boldsymbol{x}) \lambda^{d-k}(\d \boldsymbol{x})$, where $h_1 : \R^{d-k} \to \R$ is given by $h_1(x_1, \dots, x_{d-k})=2^{-(d-k)} \exp(-\sum_{j=1}^{d-k} |x_j|)$.
Choose a probability measure $\tilde \nu$ and a strictly positive measurable function $h_2:\R\to \R$ such that $\tilde \nu(\d y)=h_2(y) \nu(\d y)$. Note that $\kappa\otimes \tilde \nu (\d \boldsymbol{x}, \d y)=  h(\boldsymbol{x},y) \lambda^{d-k}(\d \boldsymbol{x}) \nu(\d y)$, where $h(\boldsymbol{x},y)=h_1(\boldsymbol{x})h_2(y)$. On some probability space $(\tilde \Omega,\tilde{\mathcal F},\tilde \P)$ let 
 $(\tilde {\boldsymbol{V}}^1_j)_{j\in \N}$ be  an i.i.d.\ sequence of $\mathcal U([0,1]^k)$-distributed random vectors, 
  $(\tilde {\boldsymbol{V}}^2_j)_{j\in \N}$ be an i.i.d.\ sequence of random vectors with the common distribution $\kappa$, 
 $(W_j)_{j\in \N}$ be  an i.i.d.\ sequence of real-valued random variables with the common distribution $\tilde \nu$,  
  and let 
 $(\Gamma_j)_{j\in \N}$ be  a sequence of partial sums of i.i.d.\ standard exponential  random variables. 
  Assume that the four sequences $(\tilde {\boldsymbol{V}}^1_j)_{j\in \N}$,   $(\tilde {\boldsymbol{V}}^2_j)_{j\in \N}$, 
 $(W_j)_{j\in \N}$ and  $(\Gamma_j)_{j\in \N}$ 
    are independent, and set 
  \begin{equation}
  \tilde J_j= W_j \1 (h(\tilde {\boldsymbol{V}}^2_j, W_j)\leq  \Gamma_j^{-1}),\  j\in \N, \qquad \text{and }\qquad \tilde \Lambda=\sum_{j=1}^\infty \delta_{(\tilde {\boldsymbol{V}}^1_j, \tilde {\boldsymbol{V}}^2_j,\tilde J_j)}.
  \end{equation}
Then $\tilde \Lambda$ is a Poisson random measure on $[0,1]^k\times \R^{d-k}\times \R_0$ 
with intensity measure $\lambda^k\otimes \lambda^{d-k}\otimes \nu$. Using again that our probability  
space is rich enough to support a  $\mathcal U([0,1])$-distributed random variable independent of $L$, we deduce by 
 Proposition~2.1 in \cite{Ros-Point} that there exists a  sequence $( {\boldsymbol{V}}^1_j,  {\boldsymbol{V}}^2_j, J_j)_{j\in \N}$
defined on $(\Omega,\mathcal F,\P)$  
which equals $(\tilde {\boldsymbol{V}}^1_j, \tilde {\boldsymbol{V}}^2_j,\tilde J_j)_{j\in \N}$  in distribution, and satisfies
\begin{equation}\label{sdljfsldhfs}
\Lambda
= \sum_{j=1}^\infty \delta_{( {\boldsymbol{V}}_j^1,  {\boldsymbol{V}}_j^2, J_j)}
\end{equation}
on $[0,1]^k \times \R^{d-k} \times \R_0$ almost surely.
In the following definition we will 
introduce the Poisson random measure appearing in the limit 
of Theorem~\ref{thm2}(i). 

\begin{definition} \label{sldjfsghsklhs} 
Let    $(\boldsymbol{U}_j)_{j\in \N}$ be an i.i.d.\ sequence of $\mathcal U([0,1]^k)$-distributed random vectors, defined on an extension of $(\Omega,\mathcal F,\P)$ and independent of $\mathcal F$, and set 
 \begin{equation}\label{ljsdfljshsdfsdfslp}
\Lambda^\ddagger =\sum_{j=1}^\infty \delta_{(\boldsymbol{U}_j,{\boldsymbol{V}}^2_j, J_j)}.
\end{equation}
\end{definition}
We note  that $\Lambda^\ddagger$ appearing in Definition~\ref{sldjfsghsklhs} is a Poisson random measure on $[0,1]^k\times \R^{d-k}\times \R_0$ with intensity measure $\lambda^k\otimes \lambda^{d-k}\otimes \nu$. Moreover,   the  Poisson random measures $\Lambda^\dagger$ and  $\Lambda^\ddagger$ appearing in Definitions~\ref{ljsdlfjsd} and \ref{sldjfsghsklhs} are neither measurable with respect to $L$ nor independent of $L$.

\subsection{Proof of Theorem \ref{thm1}(i)} \label{sec4.1}

We denote the limiting variable in Theorem \ref{thm1}(i) by $Z$. We have that $|Z| < \infty$ almost surely if 
$$
\int_{\R_0 \times (0,1)^d}\min ( 1, |y|^p H(\boldsymbol{u})) \nu (\d y) \d \boldsymbol{u} \le \int_{\R_0} \min \Big(1, |y|^p \int_{(0,1)^d} H(\boldsymbol{u}) \d \boldsymbol{u} \Big) \nu (\d y) < \infty,
$$ 
where $H(\boldsymbol{u}) := \sum_{\boldsymbol{j} \in \Z^d} |\Delta_1 h(\boldsymbol{j}-\boldsymbol{u})|^p$, $\boldsymbol{u} \in (0,1)^d$.
Indeed, $\int_{\R_0} \min(1,|y|^p) \nu (\d y)<\infty$ since $p>\beta$. Therefore, we only need to show
$\int_{(0,1)^d} H(\boldsymbol{u}) \d \boldsymbol{u} <\infty$.
For large $\| \boldsymbol{y} \|$, by rewriting $\Delta_1 h(\boldsymbol{y}) = \int_{(0,1)^d} \partial^d h(\boldsymbol{y}+\boldsymbol{v}) \d \boldsymbol{v}$ and using $|\partial^d h(\boldsymbol{y})| \le C \| \boldsymbol{y}\|^{d(\alpha-1)}$ we see that $|\Delta_1 h(\boldsymbol{y})| \le C \| \boldsymbol{y} \|^{d(\alpha-1)}$. 
By changing to spherical coordinates we know that $\int_{B_1^c (\boldsymbol{0})} \| \boldsymbol{y} \|^{d(\alpha-1)p} \d \boldsymbol{y} < \infty$ if and only if $\alpha+\frac{1}{p}<1$. 
So the integral test implies that for large $\rho>0$ there exists $C>0$ such that
$\sum_{\boldsymbol{j} \in B_\rho^c (\boldsymbol{0})} |\Delta_1 h(\boldsymbol{j}-\boldsymbol{u})|^p \le C$ for all $\boldsymbol{u} \in (0,1)^d$.
Finally, for $\| \boldsymbol{j} \|<\rho$, we have 
$\int_{(0,1)^d} |\Delta_1 h(\boldsymbol{j}-\boldsymbol{u})|^p \d \boldsymbol{u} \le C \int_{B_{2\rho}(\boldsymbol{0})} \| \boldsymbol{y} \|^{d \alpha p} \d \boldsymbol{y} < \infty$
since $\alpha+\frac{1}{p}>0$. Hence, we conclude that $|Z| < \infty$ almost surely. \\ \\
\noindent
Now, we start with the proof of the stable convergence, which is divided into two steps. In Step 1 we prove Theorem \ref{thm1}(i) for $\nu(\R_0)<\infty$, which corresponds to treatment of ``big jumps of $L$''.  In Step 2 we show that ``small jumps of $L$'' are asymptotically negligible and complete the proof of Theorem \ref{thm1}(i) for $\nu(\R_0)=\infty$. 

\bigskip
\noindent
{\it Step 1.} In the following we will prove Theorem~\ref{thm1}(i) in case where $\nu(\R_0)<\infty$. 
Choose a small $\epsilon \in (0,1)$. 
For every 
$\boldsymbol{0} \le \boldsymbol{i} < \boldsymbol{n}
$, decompose 
$\Delta_{1/n} X(\boldsymbol{i}/n) =M_{n, \epsilon} (\boldsymbol{i})+ R'_{n,\epsilon}(\boldsymbol{i})+R_{n,\epsilon}(\boldsymbol{i})$ so that 
\begin{align}
M_{n,\epsilon} (\boldsymbol{i}) &= \int_{B_\epsilon(\boldsymbol{i}/n) \cap [0,1]^d} \Delta_{1/n} g \left( \boldsymbol{i}/n - \boldsymbol{s} \right) 
L (\d \boldsymbol{s}) , \label{def:reprM}\\
R'_{n,\epsilon}(\boldsymbol{i}) &=  \int_{B_\epsilon(\boldsymbol{i}/n) \setminus [0,1]^d}
\Delta_{1/n} g \left( \boldsymbol{i}/n - \boldsymbol{s} \right) 
L (\d \boldsymbol{s}),\label{def:reprRbullet}\\
R_{n,\epsilon}(\boldsymbol{i}) &=  \int_{B_\epsilon^c(\boldsymbol{i}/n)}
\Delta_{1/n} g \left( \boldsymbol{i}/n - \boldsymbol{s} \right) 
L (\d \boldsymbol{s}).\label{def:reprR}
\end{align}
First, we will prove the stable convergence for the power variation statistics built from  $M_{n,\epsilon}(\boldsymbol{i})$ instead of the original increments $\Delta_{1/n} X(\boldsymbol{i}/n)$; later we will show that contribution of the terms $R'_{n,\epsilon}(\boldsymbol{i})$
and $R_{n,\epsilon}(\boldsymbol{i})$ is asymptotically negligible. 
Let $\Lambda$ be given by \eqref{ljsdlfjsdlj} with the representation $ \Lambda=\sum_{k=1}^\infty \delta_{(\boldsymbol{V}_k,J_k)}$ on $[0,1]^d \times \R_0$ given by \eqref{sdljflsjd}. 
 We have
\begin{align}
M_{n,\epsilon}(\boldsymbol{i}) &= \int_{[0,1]^d \times \R_0} y \Delta_{1/n} g \left( \boldsymbol{i}/n-\boldsymbol{x} \right) \1 \left( \left\| \boldsymbol{i}/n - \boldsymbol{x} \right\| < \epsilon \right) \Lambda(\d \boldsymbol{x}, \d y)\\ 
&= \sum_{k=1}^\infty  J_k \Delta_{1/n} g \left( \boldsymbol{i}/n-\boldsymbol{V}_k \right) 
\1 \left( \left\| \boldsymbol{i}/n-\boldsymbol{V}_k \right\| < \epsilon \right),  \label{sldfjlsdhfhgsl}
\end{align}
where there are at most finitely many terms in the sum in \eqref{sldfjlsdhfhgsl} 
which are different from zero, due to the fact $\nu(\R_0)<\infty$ and hence  $\Lambda([0,1]^d\times \R_0)<\infty$ almost surely.
Let us now prove that as $n \to \infty$ on the event
\begin{align*}
\Omega_\epsilon &:=  \{ \omega \in \Omega :  \|\boldsymbol{V}_{k_1} (\omega) - \boldsymbol{V}_{k_2} (\omega) \| > 2\epsilon \text{ for all }  k_1 \neq k_2 \text{ with } |J_{k_1}(\omega)|, |J_{k_2}(\omega)|\neq 0,\\
&\qquad  \qquad \text{and } \boldsymbol{V}_k (\omega) \in [\epsilon,1-\epsilon]^d  \text{ for all } k \text{ with }|J_{k}(\omega)|\neq 0 \}
\end{align*}
it holds
\begin{equation}\label{lim:tildem}
n^{d\alpha p} \sum_{\boldsymbol{0} \le \boldsymbol{i} < \boldsymbol{n}} | M_{n,\epsilon}(\boldsymbol{i}) |^p \overset{{\cal F}\textnormal{-d}}{\to} \sum_{k=1}^\infty |J_k|^p  
\sum_{\boldsymbol{j} \in \Z^d} | \Delta_{1} h( \boldsymbol{j} -\boldsymbol{U}_k ) |^p =  Z.
\end{equation}
Here $(\boldsymbol{U}_k)_{k \in \N}$ is a sequence of independent ${\cal U} ([0,1]^d)$-distributed random vectors, defined on the extension of the underlying probability space $(\Omega, {\cal F}, \P)$ and independent of the $\sigma$-algebra ${\cal F}$. 
We first note that on $\Omega_\epsilon$ every $M_{n,\epsilon}(\boldsymbol{i})$ satisfies 
either $|M_{n,\epsilon}(\boldsymbol{i})| = 0$ or $|M_{n,\epsilon}(\boldsymbol{i})| = |J_k \Delta_{1/n} g(\boldsymbol{i}/n- \boldsymbol{V}_k)|$ for some $k$.  Hence, it holds that on $\Omega_\epsilon$, 
\begin{equation}\label{eq:main}
\sum_{\boldsymbol{0}\le\boldsymbol{i}<\boldsymbol{n}} | M_{n,\epsilon}(\boldsymbol{i}) |^p = V_{n,\epsilon}, \qquad \text{where } V_{n,\epsilon} := \sum_{k=1}^\infty |J_k|^p \sum_{\boldsymbol{j} \in B_{n\epsilon} ( \{ n \boldsymbol{V}_k \} ) 
} \left| \Delta_{1/n} g \left( (\boldsymbol{j} - \{n \boldsymbol{V}_k\})/n \right) \right|^p.
\end{equation}
Since $\Omega_\epsilon \in {\cal F}$, on $\Omega_\epsilon$ the relation \eqref{lim:tildem} follows if we prove that
\begin{equation}\label{lim:Vn}
n^{d\alpha p} V_{n,\epsilon} \overset{{\cal F}\textnormal{-d}}{\to}  Z \qquad \text{as } n \to \infty.
\end{equation}
Next, we will prove for each  $k$:
\begin{equation}\label{lim:H}
n^{d \alpha p} \sum_{\boldsymbol{j} \in B_{n\epsilon} ( \{ n \boldsymbol{V}_k \} ) 
} \left| \Delta_{1/n} g \left( (\boldsymbol{j} - \{n \boldsymbol{V}_k\})/n \right) \right|^p \overset{{\cal F}\textnormal{-d}}{\to} \sum_{\boldsymbol{j} \in \Z^d} |\Delta_1 h(\boldsymbol{j}-\boldsymbol{U}_k)|^p = H(\boldsymbol{U}_k).
\end{equation}
Under Assumption (H1) we have the identity 
$$n^{d\alpha} g \left( (\boldsymbol{j} - \{ n \boldsymbol{V}_k \})/n \right) = h(\boldsymbol{j} - \{ n \boldsymbol{V}_k \}) f 
\left( (\boldsymbol{j} - \{ n \boldsymbol{V}_k \})/n \right)
$$
with $\lim_{\boldsymbol{x} \to \boldsymbol{0}} f(\boldsymbol{x}) =1$ and by \eqref{braconv}
$$
\{ n \boldsymbol{V}_k \} \overset{{\cal F}\textnormal{-d}}{\to} \boldsymbol{U}_k.
$$
By the continuous mapping theorem for stable convergence, we get that 
\begin{align*}
n^{d \alpha p} \sum_{\boldsymbol{j} \in B_r(\boldsymbol{0})} \left|\Delta_{1/n} g\left((\boldsymbol{j}-\{n\boldsymbol{V}_k\})/n \right) \right|^p \overset{{\cal F}\textnormal{-d}}{\to} \sum_{\boldsymbol{j} \in B_r(\boldsymbol{0})} |\Delta_1 h(\boldsymbol{j}-\boldsymbol{U}_k)|^p =: H_r (\boldsymbol{U}_k)
\end{align*}
for some large $r>0$. 
Since $\lim_{r \to \infty} H_r(\boldsymbol{u}) = H(\boldsymbol{u})$ for $\boldsymbol{u} \in (0,1)^d$, it suffices to show that
\begin{equation}\label{lim:Hrem}
\lim_{r \to \infty} \limsup_{n \to \infty} \sup_{\boldsymbol{u} \in (0,1)^d} n^{d \alpha p} \sum_{\boldsymbol{j} \in B_{n\epsilon}(\boldsymbol{u}) \setminus B_r(\boldsymbol{0})} |\Delta_{1/n} g((\boldsymbol{j}-\boldsymbol{u})/n)|^p = 0.
\end{equation} 
Indeed, for 
$\boldsymbol{j} \in B_{n \epsilon} (\boldsymbol{u}) \setminus B_r (\boldsymbol{0})$, rewriting
$n^{d\alpha} \Delta_{1/n} g ( (\boldsymbol{j}-\boldsymbol{u})/n )=n^{d(\alpha-1)}  \int_{(0,1)^d} \partial^d g ( (\boldsymbol{j} - \boldsymbol{u} + \boldsymbol{v})/n ) \d \boldsymbol{v}$ with
$n^{d(\alpha-1)} |\partial^d g(( \boldsymbol{j}-\boldsymbol{u}+\boldsymbol{v} )/n)|\le C \| \boldsymbol{j}-\boldsymbol{u}+\boldsymbol{v} \|^{d (\alpha-1)}$, we get
\begin{equation}\label{ineq:Deltag1}
n^{d\alpha} \left| \Delta_{1/n} g  ( (\boldsymbol{j}-\boldsymbol{u})/n ) \right| \le C \| \boldsymbol{j} \|^{d(\alpha-1)}.
\end{equation}
Finally, we have
$\lim_{r \to \infty} \sum_{\boldsymbol{j} \in B^c_r(\boldsymbol{0}) 
} \| \boldsymbol{j} \|^{d(\alpha-1)p} = 0$ since $\alpha + 1/p < 1$, which implies \eqref{lim:Hrem} and thus completes the proof of \eqref{lim:H}. By independence and the continuous mapping theorem we get for  all $K=1,2,\dots$ 
\begin{equation}\label{sdfljsldhs}
 \sum_{k=1}^K |J_k|^p \sum_{\boldsymbol{j} \in B_{n\epsilon} ( \{ n \boldsymbol{V}_k \} ) 
} \left| \Delta_{1/n} g \left( (\boldsymbol{j} - \{n \boldsymbol{V}_k\})
/n \right) \right|^p \overset{{\cal F}\textnormal{-d}}{\to}
\sum_{k=1}^K |J_k|^p  
H(\boldsymbol{U}_k).
\end{equation}
Since the event  $A_K:=\{\omega\in \Omega: J_k(\omega)=0\text{ for all }k>K\}$ is $\mathcal F$-measurable, it follows by \eqref{sdfljsldhs}  that  \eqref{lim:Vn} holds on $A_K$ for all $K=1,2,\dots$, and since $A_K\uparrow \Omega$ as $K\to \infty$ we deduce that \eqref{lim:Vn} holds.

Next, let us prove that the terms $R'_{n,\epsilon}(\boldsymbol{i})$ in \eqref{def:reprRbullet} satisfy 
\begin{equation}\label{lim:Mprime}
\lim_{\epsilon \downarrow 0} \limsup_{n \to \infty} \P \Big( n^{d \alpha p} \sum_{\boldsymbol{0} \le \boldsymbol{i} < \boldsymbol{n}} | R'_{n,\epsilon}(\boldsymbol{i}) |^p > \delta \Big) =0
\end{equation}
for all $\delta >0$. 
For this purpose, choose a large rectangle $B'$ in $\R^d$. 
 Recall $\Lambda$ associated to $L$ by \eqref{ljsdlfjsdlj} and use its representation $\Lambda=\sum_{k=1}^\infty \delta_{(\boldsymbol{V}'_k,J'_k)}$ on $B' \setminus [0,1]^d \times \R_0$, analogous to that in \eqref{sdljflsjd}. Then, for $\bar p=\max(p,1)$, it holds that
\begin{align}\label{ineq:calM'}
\Big( n^{d \alpha p} \sum_{\boldsymbol{0} \le \boldsymbol{i} < \boldsymbol{n}} | R'_{n,\epsilon}(\boldsymbol{i}) |^p \Big)^{1/\bar p} \le 
\sum_{k=1}^{\infty} \big( G_{n,\epsilon} ( \boldsymbol{V}'_k ) | J'_k |^{p} \big)^{1/\bar p} 
\end{align}
with
$$
G_{n,\epsilon} (\boldsymbol{V}'_k) := n^{d \alpha p} \sum_{\boldsymbol{0} \le \boldsymbol{i} < \boldsymbol{n}} \left| \Delta_{1/n} g ( \boldsymbol{i}/n - \boldsymbol{V}'_k ) \right|^p \1 \left( \| \boldsymbol{i}/n - \boldsymbol{V}'_k \| < \epsilon \right).
$$ 
Note that a ${\cal U}(B' \setminus [0,1]^d)$-distributed random vector $\boldsymbol{V}'_k$ does not belong to $B_\epsilon(\boldsymbol{i}/n)$ if $\boldsymbol{i} \in [n\epsilon, n(1-\epsilon)]^d$. Therefore, $\E [|G_{n,\epsilon}(\boldsymbol{V}'_k)|] \le C \epsilon (I_{n,\epsilon}^0+I_{n,\epsilon}^1)$
with 
\begin{equation}\label{def:I01}
I^0_{n,\epsilon} := n^{pd \alpha + d} \int_{\| \boldsymbol{x} \| < \frac{d}{n}} | \Delta_{\frac{1}{n}} g(\boldsymbol{x})|^p \d \boldsymbol{x}, \quad
I^1_{n,\epsilon} := n^{pd \alpha + d} \int_{\frac{d}{n} \le \| \boldsymbol{x} \| < \epsilon} | \Delta_{\frac{1}{n}} g(\boldsymbol{x}) |^p \d \boldsymbol{x},
\end{equation}
where
\begin{align*}
I^0_{n,\epsilon} \le C n^{pd \alpha+d} \int_{\| \boldsymbol{x} \|
	<\frac{2d}{n}} |g(\boldsymbol{x})|^p \d \boldsymbol{x} \le C n^{p d \alpha +d} \int_{\| \boldsymbol{x} \|<\frac{2d}{n}} \| \boldsymbol{x} \|^{pd \alpha} \d \boldsymbol{x} \le C \int_{\| \boldsymbol{x} \| < 2d} \| \boldsymbol{x} \|^{pd \alpha} \d \boldsymbol{x} < \infty
\end{align*}
since $pd \alpha +d -1>-1$, whereas $|n^{d} \Delta_{\frac{1}{n}} g(\boldsymbol{x})|  \le \int_{(0,1)^d} | \partial^d g(\boldsymbol{x}+\frac{\boldsymbol{u}}{n}) | \d \boldsymbol{u}
\le C \| \boldsymbol{x} \|^{d(\alpha-1)}$ for $\frac{d}{n} \le \| \boldsymbol{x} \| < \epsilon$. This implies  
\begin{align*}
I^1_{n,\epsilon} &\le C n^{p d\alpha + d - pd} \int_{\frac{d}{n} \le \| \boldsymbol{x} \| < \epsilon } \| \boldsymbol{x} \|^{pd(\alpha-1)} \d \boldsymbol{x} \le C \int_{d \le \| \boldsymbol{x} \|} \| \boldsymbol{x} \|^{pd(\alpha-1)} \d \boldsymbol{x} < \infty
\end{align*}
since $pd \alpha - pd + d - 1 < -1$. 
From $\E[|G_{n,\epsilon}(\boldsymbol{V}'_k)|] \le C \epsilon$ it follows 
$$
\lim_{\epsilon\downarrow0} \limsup_{n \to \infty} \P \big( G_{n,\epsilon} ( \boldsymbol{V}'_k ) | J'_k|^p > \delta \big) = 0,
$$
hence
$$
\lim_{\epsilon \downarrow 0} \limsup_{n \to \infty} \P \Big( \sum_{k=1}^\infty \big( G_{n,\epsilon} ( \boldsymbol{V}'_k ) | J'_k |^p \big)^{1/\bar p} > \delta^{1/\bar p} \Big) = 0,
$$
which in turn implies \eqref{lim:Mprime}.

Finally, we consider the terms $R_{n,\epsilon}(\boldsymbol{i})$ having representation \eqref{def:reprR}. We prove that  
\begin{equation}\label{lim:R}
n^{pd \alpha} \sum_{\boldsymbol{0} \le \boldsymbol{i} < \boldsymbol{n}} |R_{n,\epsilon}(\boldsymbol{i})|^p \overset{\P}{\to} 0 
\qquad \text{as } n\to \infty.
\end{equation}
For this purpose, we will first determine a bounded function $\psi \in L^\theta (\R^d)$, which satisfies
\begin{equation}\label{ineq:psi}
n^d | \Delta_{1/n} g(\boldsymbol{i}/n-\boldsymbol{x})| 
\1 (\boldsymbol{x} \in B_\epsilon^c (\boldsymbol{i}/n)) 
\le \psi (\boldsymbol{x})
\end{equation}
for all $\boldsymbol{x} \in \R^d$, $\boldsymbol{0} \le \boldsymbol{i}<\boldsymbol{n}$ 
and large enough $n \in \N$. Let $\rho >0$ be large. Consider the identity
$$
n^d \Delta_{1/n} g( \boldsymbol{i}/n - \boldsymbol{x} ) = \int_{(0,1)^d} \partial^d g( (\boldsymbol{u} + \boldsymbol{i})/n - \boldsymbol{x} ) \d \boldsymbol{u},
$$ 
where $|\partial^d g((\boldsymbol{u}+\boldsymbol{i})/n - \boldsymbol{x})| \le C (\epsilon/2)^{d(\alpha-1)}$ if $\boldsymbol{x} \in B_{2\rho}(\boldsymbol{0}) \cap B_\epsilon^c (\boldsymbol{i}/n)$, since $|\partial^d g ( \boldsymbol{v} )| \le C \| \boldsymbol{v} \|^{d(\alpha-1)}$, $\boldsymbol{v} \in B_{3\rho} (\boldsymbol{0})$. Furthermore,  $|\partial^d g((\boldsymbol{u}+\boldsymbol{i})/n - \boldsymbol{x})| \le |\partial^d g(\boldsymbol{x}/2)|$ if $\boldsymbol{x} \in B_{2\rho}^c (\boldsymbol{0})$, by monotonicity of $|\partial^d g|$ on $B_\rho^c (\boldsymbol{0})$.
Consequently, for $\boldsymbol{x} \in \R^d$, we define
$$
\psi (\boldsymbol{x}) := C 
\1 ( \boldsymbol{x} \in B_{2\rho} (\boldsymbol{0}) ) 
+ | \partial^d g(\boldsymbol{x}/2) | 
\1 ( \boldsymbol{x} \in B_{2\rho}^c (\boldsymbol{0}) ),
$$
where $C$ depends on $\epsilon$. In what follows, w.l.o.g.\ assume $|\psi (\boldsymbol{x}) | \le 1$, $\boldsymbol{x} \in \R^d$. 

With $\Lambda$ given by \eqref{ljsdlfjsdlj} we set $\Lambda^1 (\cdot)=\Lambda(\cdot\cap \{ (\boldsymbol{x},y) \in \R^d \times \R_0 : |\psi(\boldsymbol{x})y| >  1\})$ and for all $B\in {\cal B}_b(\R^d)$ set 
$$
L^1 (B) = \int_{B \times\R_0} y \Lambda^1 (\d \boldsymbol{x}, \d y)  \qquad \text{and} \qquad L^0 (B) = L(B)- L^1(B).
$$
The $L^0$ and $L^1$ are independent infinitely divisible random measures such that for every $B \in {\cal B}_b (\R^d)$,
\begin{align*}
\E \left[ \e^{\i t L^0 (B)} \right] &= \exp \ \Big( \int_{B \times \R_0} ( \e^{\i t y} - 1 - \i t y\1(|y|\le 1) ) \1 (|\psi(\boldsymbol{x})y| \le 1) \d \boldsymbol{x} \nu (\d y) \Big),\\
\E \Big[ \e^{\i t L^1 (B)} \Big] &=\exp \Big( \int_{B \times \R_0} ( \e^{\i t y} - 1 - \i t y\1(|y|\le 1) ) \1 (|\psi(\boldsymbol{x})y| > 1) \d \boldsymbol{x} \nu (\d y) \Big), \qquad t \in \R.
\end{align*}
Then, for every $\boldsymbol{i}$, we decompose $n^d R_{n,\epsilon}(\boldsymbol{i}) = Q^0_{n,\epsilon}(\boldsymbol{i}) + Q^1_{n,\epsilon}(\boldsymbol{i})$, where
\begin{align*}
Q^j_{n,\epsilon}(\boldsymbol{i}) := &:= \int_{B^c_\epsilon(\boldsymbol{i}/n)} n^d \Delta_{1/n} g( \boldsymbol{i}/n - \boldsymbol{s} ) L^j (\d \boldsymbol{s}), \qquad j=0,1.
\end{align*}
We claim that for $j=0,1$, 
$$
n^{pd(\alpha - 1)} \sum_{\boldsymbol{0} \le \boldsymbol{i} < \boldsymbol{n}} |Q^j_{n,\epsilon} (\boldsymbol{i}) |^p \overset{\P}{\to} 0
$$
follows from
\begin{equation}\label{bdd:Q}
\sup_{n \in \N, \, \boldsymbol{0} \le \boldsymbol{i} < \boldsymbol{n}} \E [|Q_{n,\epsilon}^0 (\boldsymbol{i})|^p] < \infty \qquad \text{and} \qquad \sup_{n \in \N, \, \boldsymbol{0} \le \boldsymbol{i} < \boldsymbol{n}} |Q_{n,\epsilon}^1 (\boldsymbol{i})| < \infty \text{ a.s.},
\end{equation} 
since $pd\alpha -pd+d <0$. 

For the first relation in \eqref{bdd:Q}, it suffices to show that 
$$
\int_{B^c(\boldsymbol{i}/n)} \Phi_p \big(|n^d \Delta_{1/n} g(\boldsymbol{i}/n - \boldsymbol{x})|, \boldsymbol{x} \big) \d \boldsymbol{x} \le C,
$$
where 
$$
\Phi_p (v,\boldsymbol{x}) = \int_{\R_0} \big( |vy|^p\1 (|vy| >1) + |vy|^2 \1 (|vy| \le 1) \big) \1 (|\psi(\boldsymbol{x})y| \le 1) \nu (\d y),
$$
cf.\ Theorem 3.3 in \cite{rajput1989}.
In view of \eqref{ineq:psi} we have that
\begin{align*}
\int_{B^c (\boldsymbol{i}/n)} \Phi_p \big( |n^d \Delta_{1/n} g(\boldsymbol{i}/n - \boldsymbol{x}) |, \boldsymbol{x} \big) \d \boldsymbol{x} \le \int_{\R^d \times \R_0} | \psi(\boldsymbol{x}) y |^2 \1 (|\psi(\boldsymbol{x})y| \le 1) \nu (\d y),
\end{align*}
where the estimate \eqref{ineq:estV} implies
\begin{align*}
\int_{\R_0} |xy|^2 \1 (|xy| \le 1) \nu (\d y) \le C |x|^\theta \qquad \text{for } |x| \le 1,
\end{align*}
and $\psi \in L^\theta (\R^d)$ is bounded. We conclude that the first relation in \eqref{bdd:Q} holds.
Finally, the second relation in \eqref{bdd:Q} follows in view of \eqref{ineq:psi} from
\begin{align*}
|Q_{n,\epsilon}^1 (\boldsymbol{i})| \le \int_{B^c_\epsilon (\boldsymbol{i}/n) \times \R_0} |n^d \Delta_{1/n} g (\boldsymbol{i}/n - \boldsymbol{x}) y|  \Lambda^1 (\d \boldsymbol{x}, \d y) \le \int_{\R^d \times \R_0} |\psi (\boldsymbol{x}) y | \Lambda^1 (\d \boldsymbol{x}, \d y) < \infty,
\end{align*}
where the last stochastic integral is well-defined because we have that $\psi \in L^\theta (\R^d)$ is bounded and
$$
\int_{\R_0} \min (|xy|,1) \1(|x y| >1) \nu (\d y) = 2 \int_0^\infty \1 (|xy|>1) \nu (\d y) \le C |x|^\theta \qquad \text{for } |x| \le 1
$$
by \eqref{ineq:estV}. This completes the proof of \eqref{lim:R}.

Let us now complete the proof of Theorem \ref{thm1}(i) in case $\nu(\R_0)<\infty$. For some small $\epsilon \in (0,1)$ we have the decomposition  $\Delta_{\frac{1}{n}}X(\frac{\boldsymbol{i}}{n}) = M_{n,\epsilon}(\boldsymbol{i})+R'_{n,\epsilon}(\boldsymbol{i})+R_{n,\epsilon}(\boldsymbol{i})$. Correspondingly, with $\bar p := \max(p,1)$ we decompose
\begin{align*}
(V_n(p) )^{\frac{1}{\bar p}} &= ( V_n(p) )^{\frac{1}{\bar p}} - \Big( \sum_{\boldsymbol{0} \le \boldsymbol{i} < \boldsymbol{n} } |M_{n,\epsilon}(\boldsymbol{i})|^p \Big)^{\frac{1}{\bar p}} + \Big( \sum_{\boldsymbol{0}\le \boldsymbol{i}<\boldsymbol{n}} |M_{n,\epsilon}(\boldsymbol{i})|^p \Big)^{\frac{1}{\bar p}}.
\end{align*}
Concerning the last term, the limiting result \eqref{lim:tildem} holds on the event $\Omega_\epsilon$ with the limit satisfying $Z^{\frac{1}{\bar p}} \1(\Omega_\epsilon) \to Z^{\frac{1}{\bar p}}$, since $\P(\Omega_\epsilon) \uparrow 1$ as $\epsilon \downarrow 0$.
Applying \eqref{ineq0} and \eqref{ineq1}, we see that
\begin{align*}
\Big| ( V_n(p) )^{\frac{1}{\bar p}} - \Big( \sum_{\boldsymbol{0} \le \boldsymbol{i} < \boldsymbol{n} } |M_{n,\epsilon}(\boldsymbol{i})|^p \Big)^{\frac{1}{\bar p}} \Big|
&\le\Big( \sum_{\boldsymbol{0} \le \boldsymbol{i}<\boldsymbol{n}} |R'_{n,\epsilon}(\boldsymbol{i})|^p \Big)^{\frac{1}{\bar p}} + \Big( \sum_{\boldsymbol{0} \le \boldsymbol{i}<\boldsymbol{n}} |R_{n,\epsilon}(\boldsymbol{i})|^p \Big)^{\frac{1}{\bar p}}.
\end{align*}
where the r.h.s.\ terms satisfy \eqref{lim:Mprime}, \eqref{lim:R}, proving that 
$$
\lim_{\epsilon \downarrow 0} \limsup_{n \to \infty} \P\Big( \Big| (n^{d \alpha p} V_n(p) )^{\frac{1}{\bar p}} - \Big( n^{d \alpha p} \sum_{\boldsymbol{0} \le \boldsymbol{i} < \boldsymbol{n} } |M_{n,\epsilon}(\boldsymbol{i})|^p \Big)^{\frac{1}{\bar p}} \Big| > \delta \Big) =0
$$
for all $\delta>0$. We conclude that $(n^{d \alpha p} V_n(p))^{\frac{1}{\bar p}} \overset{{\cal F}\textnormal{-d}}{\to} Z^{\frac{1}{\bar p}}$ as $n \to \infty$. 
\bigskip

\noindent
{\it Step 2.} Let $\nu (\R_0) = \infty$. We choose a some small $\epsilon>0$, 
and use $\Lambda$ given by \eqref{ljsdlfjsdlj} to define $\Lambda^{> \epsilon}(\cdot)=\Lambda(\cdot\cap (\R^d\times [-\epsilon,\epsilon]^c))$ and for all $B\in {\cal B}_b(\R^d)$ set 
\begin{align}
L^{> \epsilon} (B)=   \int_{B\times \R_0} y\Lambda^{> \epsilon}(\d \boldsymbol{x}, \d y)\qquad \text{and} \qquad L^{\le \epsilon}(B)=L(B)-L^{> \epsilon}(B). 
\end{align}
Then $L^{\le \epsilon}$ and $L^{>\epsilon}$ are  independent infinitely divisible random measures such that for every $B \in {\cal B}_b(\R^d)$,
\begin{align*}
\E\left[ \e^{\i t L^{\le \epsilon} (B)}\right] &= \exp \Big( \lambda^d (B) \int_{0<|y|\le \epsilon} (\e^{\i t y} - 1 - \i t y \1 (|y| \le 1)) \nu (\d y) \Big),\\
\E\left[ \e^{\i t L^{>\epsilon}(B)}\right] &= \exp \Big( \lambda^d (B) \int_{|y|>\epsilon} (\e^{\i t y} - 1 - \i t y \1 (|y| \le 1)) \nu (\d y) \Big), \qquad t \in \R.
\end{align*}
Then we decompose  $\Delta_{1/n} X(\boldsymbol{i}/n)=\Delta_{1/n} X^{\le \epsilon}(\boldsymbol{i}/n)+\Delta_{1/n} X^{>\epsilon} (\boldsymbol{i}/n)$ with
$$
\Delta_{1/n} X^{\le \epsilon}(\boldsymbol{i}/n) =\int_{\R^d} \Delta_{1/n} g ( \boldsymbol{i}/n - \boldsymbol{s} ) L^{\le \epsilon} (\d \boldsymbol{s})  \qquad \text{and} \qquad 
\Delta_{1/n} X^{>\epsilon}(\boldsymbol{i}/n)=\int_{\R^d} \Delta_{1/n} g ( \boldsymbol{i}/n - \boldsymbol{s} ) L^{>\epsilon} (\d \boldsymbol{s}).
$$
Let $\Lambda^\dagger$ be the Poisson random measure given by \eqref{ljsdfljshslp}.  
Since $\nu( [-\epsilon,\epsilon]^c ) < \infty$, we obtain by Step~1 that 
$$
n^{d \alpha p} \sum_{\boldsymbol{0} \le \boldsymbol{i} < \boldsymbol{n}} |\Delta_{1/n} X^{>\epsilon} ( \boldsymbol{i}/n ) |^p \overset{{\cal F}\textnormal{-d}}{\to} \int_{[0,1]^d \times [-\epsilon,\epsilon]^c} |y|^p \sum_{\boldsymbol{j} \in \Z^d} |\Delta_{1}h(\boldsymbol{j}-\boldsymbol{u})|^p \Lambda^\dagger (\d \boldsymbol{u}, \d y) =: Z^{>\epsilon}  \qquad \text{ as } n \to \infty.
$$
On the other hand,  as $\epsilon \downarrow 0$
$$
Z^{>\epsilon} \overset{\P}{\to} \int_{[0,1]^d \times \R_0} |y|^p \sum_{\boldsymbol{j} \in \Z^d} |\Delta_{1}h(\boldsymbol{j}-\boldsymbol{u})|^p \Lambda^\dagger (\d \boldsymbol{u}, \d y) = Z.
$$
By \eqref{ineq0} and \eqref{ineq1}, it only remains to show that for all $\delta >0$,
\begin{equation}\label{lim:P}
\lim_{\epsilon \downarrow 0} \limsup_{n \to \infty} \P \Big( n^{d\alpha p} \sum_{\boldsymbol{0} \le \boldsymbol{i} < \boldsymbol{n}} | \Delta_{1/n} X^{\le \epsilon} ( \boldsymbol{i}/n) |^p > \delta \Big) = 0.
\end{equation}
Indeed, by  Markov's inequality \eqref{lim:P} follows  if we  show that 
$$
\lim_{\epsilon \downarrow 0} \limsup_{n \to \infty} n^{d \alpha p + d} \E[ |\Delta_{\frac{1}{n}} X^{\le \epsilon} (\boldsymbol{0})|^p] = 0,
$$
for which it suffices to show that
\begin{equation}\label{lim:P2}
\lim_{\epsilon \downarrow 0} \limsup_{n \to \infty} \int_{\R^d} \int_{0<|y|\le \epsilon} \phi_p ( n^{d(\alpha+\frac{1}{p})} \Delta_{\frac{1}{n}} g(\boldsymbol{x}) y ) \nu (\d y) \d \boldsymbol{x} = 0,
\end{equation}
where 
$\phi_p(y) := |y|^p \1 (|y|>1) + |y|^2 \1(|y| \le 1)$ for $y \in \R$, cf.\ Theorem 3.3 in \cite{rajput1989}. 
Using \eqref{ineq:psi} with bounded $\psi \in L^\theta(\R^d)$, we obtain
\begin{align*}
&\int_{B^c_1(\boldsymbol{0})} \int_{0<|y|\le\epsilon} \phi_p (n^{d(\alpha+\frac{1}{p})} \Delta_{\frac{1}{n}} g(\boldsymbol{x}) y) \nu (\d y) \d \boldsymbol{x}\\ 
&\qquad \le \int_{B^c_1(\boldsymbol{0})} \int_{0<|y|\le\epsilon} \Big( |n^{d(\alpha+\frac{1}{p}-1)} \psi(\boldsymbol{x}) y|^p\1(|\psi(\boldsymbol{x})|>1) + | n^{d(\alpha+\frac{1}{p}-1)} \psi(\boldsymbol{x}) y|^2 \Big) \nu (\d y) \d \boldsymbol{x} = o(1)
\end{align*}
as $n \to \infty$, since $\alpha+\frac{1}{p}<1$. Using $\phi_p(y) \le |y|^p + |y|^2 \1 (p>2)$ for $y \in \R$, we get
\begin{align*}
\int_{B_1(\boldsymbol{0})} \int_{0<|y|\le\epsilon} \phi_p ( n^{d(\alpha+\frac{1}{p})} \Delta_{\frac{1}{n}} g(\boldsymbol{x}) y ) \nu (\d y)  \d \boldsymbol{x} \le 
\displaystyle I_n (p) \int_{0<|y|\le\epsilon} |y|^p \nu (\d y) + I_n (2) \int_{0<|y|\le\epsilon} |y|^2 \nu (\d y) \1 (p >2)
\end{align*}
with the second term present on the r.h.s.\ only if $p>2$ and with
$$
I_n (q) := \int_{B_1(\boldsymbol{0})} | n^{d (\alpha + \frac{1}{p})} \Delta_{\frac{1}{n}} g(\boldsymbol{x}) |^q \d \boldsymbol{x}, \qquad q>0.
$$
Note that by Jensen's inequality $I_n (2) \le C (I_n(p))^{\frac{2}{p}}$ if $p>2$, whereas $I_n (p) \le C$ follows from analysis of the integrals in \eqref{def:I01}. Similarly to \eqref{ineq:numom0}, we have $\int_0^{\epsilon} y^p \nu(\d y) \le C \epsilon^{p-\beta} = o(1)$ as $\epsilon \downarrow 0$, since $p>\beta$. This completes the proof of \eqref{lim:P2} and \eqref{lim:P}, and hence the proof of Theorem \ref{thm1}(i).

\subsection{Proof of Theorem \ref{thm1}(ii)}

Let us  verify that the limiting constant in Theorem \ref{thm1}(ii)
is finite. This follows from 
\begin{equation}\label{lim1exists}
\int_{\R^d} |\Delta_1 h (\boldsymbol{s})|^\beta \d s < \infty.
\end{equation}
Choose $\rho>0$ to be large. 
Then it holds
$$
\int_{B_{2\rho}(\boldsymbol{0})} |\Delta_1 h(\boldsymbol{s})|^\beta \le C \int_{B_{3\rho} (\boldsymbol{0})} |h(\boldsymbol{s})|^\beta \d \boldsymbol{s} = C \int_{B_{3\rho} (\boldsymbol{0})} \| \boldsymbol{s} \|^{d \alpha \beta} \d \boldsymbol{s}
< \infty
$$
if and only if $\alpha>-1/\beta$.
For $\boldsymbol{s} \in B^c_{2\rho}(\boldsymbol{0})$, rewrite
$$
\Delta_1 h(\boldsymbol{s}) = \int_{[0,1)^d} \partial^d h (\boldsymbol{s}+\boldsymbol{u}) \d \boldsymbol{u},
$$
where $\partial^d h (\boldsymbol{s}) = \| \boldsymbol{s} \|^{d(\alpha-1)} \ell (\boldsymbol{s})$ with $\ell (\boldsymbol{s}) := \prod_{i=1}^{d} (d\alpha-2(i-1)) (s_i/\| \boldsymbol{s} \|)$ satisfies
$$
|\Delta_1 h (\boldsymbol{s})| \le \int_{[0,1)^d} |\partial^d h(\boldsymbol{s} + \boldsymbol{u})| \d \boldsymbol{u} \le  C \int_{[0,1)^d} \| \boldsymbol{s} + \boldsymbol{u} \|^{d(\alpha-1)} \d \boldsymbol{u} \le C \| \boldsymbol{s} \|^{d (\alpha-1)}.
$$
Then
$$
\int_{B_{2\rho}^c(\boldsymbol{0})} \| \boldsymbol{s} \|^{d(\alpha-1)\beta} \d \boldsymbol{s}
< \infty
$$
if and only if $\alpha+1/\beta<1$. Hence, \eqref{lim1exists} holds.

Now, we show the convergence in probability in Theorem \ref{thm1}(ii). 
Using the scaling property of the $\beta$-stable random measure $L$, we have that 
$\{ n^{dH} \Delta_{1/n} X(\boldsymbol{i}/n) \}_{\boldsymbol{i} \in \Z^d} \overset{{\rm fdd}}{=} \{ Y_n (\boldsymbol{i})\}_{\boldsymbol{i} \in \Z^d }$
with 
$$
Y_n (\boldsymbol{i}) := \int_{\R^d} n^{d\alpha} \Delta_{1/n} g((\boldsymbol{i} - \boldsymbol{s})/n) L (\d \boldsymbol{s}).
$$
Thus, we deduce the distributional identity
\begin{equation}\label{eq:V}
n^{dHp} V_n (p) \overset{{\rm d}}{=} \sum_{\boldsymbol{0} \le \boldsymbol{i} < \boldsymbol{n}} |Y_n(\boldsymbol{i})|^p.
\end{equation}
Next, we approximate $(Y_n (\boldsymbol{i}) )_{\boldsymbol{i} \in \Z^d}$ by $( Y_\infty (\boldsymbol{i}) )_{\boldsymbol{i} \in \Z^d}$, where
$$
Y_\infty (\boldsymbol{i}) := \int_{\R^d} \Delta_1 h ( \boldsymbol{i} - \boldsymbol{s} ) L( \d \boldsymbol{s} )
$$
is well defined due to \eqref{lim1exists}. More specifically, we will prove that 
\begin{equation}\label{lim:Yn-Yinfty}
\E [| Y_n (\boldsymbol{0}) - Y_\infty (\boldsymbol{0}) |^p] = C \Big( \int_{\R^d} | n^{d\alpha} \Delta_{1/n} g ( \boldsymbol{s}/n ) - \Delta_1 h ( \boldsymbol{s} )  |^\beta \d \boldsymbol{s} \Big)^{p/\beta} =o(1).
\end{equation}
Observe that for almost every $\boldsymbol{s} \in \R^d$, the pointwise convergence $n^{d\alpha} \Delta_{1/n} g ( \boldsymbol{s}/n ) \to \Delta_1 h ( \boldsymbol{s} )$ follows from the definition of $g$ and homogeneity of $h$. Let us verify the dominated convergence argument. By the definition of $g$ and homogeneity of $h$, we have $n^{d\alpha} |g(\boldsymbol{s}/n)| \le C \max (1, \| \boldsymbol{s}\|^{d \alpha})$ for $\| \boldsymbol{s} \| < 3\rho$. For $2\rho \le \| \boldsymbol{s} \|< 2\rho n$, we have
$$
n^{d\alpha} |\Delta_{1/n} g ( \boldsymbol{s}/n )| \le n^{d(\alpha-1)} \int_{[0,1)^d} | \partial^d g ( (\boldsymbol{s}+\boldsymbol{u})/n ) | \d \boldsymbol{u} \le C \| \boldsymbol{s} \|^{d(\alpha-1)}
$$
using $\| \boldsymbol{s} + \boldsymbol{u} \| \ge \| \boldsymbol{s} \|/2$ and $|\partial^d g(\boldsymbol{v})| \le C \| \boldsymbol{v} \|^{d(\alpha-1)}$, $\| \boldsymbol{v} \| < 3\rho$. Hence the dominated convergence theorem in $L^\beta (\R^d)$ implies
$$
\int_{B_{2\rho n}(\boldsymbol{0})} | n^{d\alpha} \Delta_{1/n} g ( \boldsymbol{s}/n ) - \Delta_1 h ( \boldsymbol{s} )  |^\beta \d \boldsymbol{s} = o (1).
$$
We next consider
\begin{align*}
I_{n} &:= n^{d\alpha\beta} \int_{B^c_{2\rho n}(\boldsymbol{0})} | \Delta_{1/n} g(\boldsymbol{s}/n) |^\beta \d \boldsymbol{s},
\end{align*}
where
$$
n^{d} | \Delta_{1/n} g(\boldsymbol{s}/n) | \le \int_{[0,1)^d} | \partial^d g ((\boldsymbol{s} + \boldsymbol{u}) / n) | \d \boldsymbol{u} \le | \partial^d g (\boldsymbol{s}/(2n)) | 
$$
using $\| (\boldsymbol{s} + \boldsymbol{u})/n \| \ge  \| \boldsymbol{s}/(2n) \| \ge \rho$ and
the monotonicity of $|\partial^d g|$ on $B^c_{\rho} (\boldsymbol{0})$. Hence
\begin{align*}
I_{n} &\le C n^{d(\alpha -1)\beta} \int_{\R^d} | \partial^d g ( \boldsymbol{s} / (2n) ) |^\beta \1 ( \| \boldsymbol{s} \| \ge 2\rho n ) \d \boldsymbol{s}\\
&= C n^{d (H - 1) \beta} \int_{B^c_\rho (\boldsymbol{0})} | \partial^d g(\boldsymbol{s}) |^\beta \d \boldsymbol{s} =o(1),
\end{align*}
since $H < 1$. From this estimate and \eqref{lim1exists} it follows that
$$
\int_{B^c_{2\rho n} (\boldsymbol{0})} | n^{d\alpha} \Delta_{1/n} g ( \boldsymbol{s}/n ) - \Delta_1 h ( \boldsymbol{s} )  |^\beta \d \boldsymbol{s} = o (1).
$$
This completes the proof of \eqref{lim:Yn-Yinfty}, which implies convergence in probability
\begin{equation}\label{lim:approx}
n^{-d} \sum_{\boldsymbol{0} \le \boldsymbol{i} < \boldsymbol{n}} | Y_n (\boldsymbol{i}) - Y_\infty (\boldsymbol{i}) |^p \overset{\P}{\to} 0.
\end{equation}
By combining Theorem~4.1 and Remark~4.3 of \cite{WRS} it follows that the stationary process $(Y_\infty(\boldsymbol{i}))_{\boldsymbol{i}\in \Z^d}$
is ergodic since it is a  stable moving average. Therefore,  we obtain from   a multiparameter Birkhoff theorem \cite[Theorem 2.8]{WRS} the convergence 
\begin{equation}\label{lim:ergodic}
n^{-d} \sum_{\boldsymbol{0} \le \boldsymbol{i} < \boldsymbol{n} } | Y_\infty (\boldsymbol{i}) |^p  \overset{\P}{\to} \E[ | Y_\infty (\boldsymbol{0}) |^p].
\end{equation}
By \eqref{ineq0}, \eqref{ineq1}, \eqref{lim:approx} and \eqref{lim:ergodic} it follows that 
$$
n^{d(Hp-1)} V_n(p) \overset{{\rm d}}{=} n^{-d} \sum_{\boldsymbol{0} \le \boldsymbol{i} < \boldsymbol{n}} |Y_n (\boldsymbol{i})|^p  
\overset{\P}{\to} \E [|Y_\infty (\boldsymbol{0})|^p].
$$
Due to the  scaling properties of stable random variables it follows that $\E [| Y_\infty (\boldsymbol{0}) |^p]$  
coincides with the limiting constant in the statement of Theorem \ref{thm1}(ii), and hence the proof of convergence in probability is complete.

Finally, we recall that convergence in $L^1$ follows from convergence in probability and uniform integrability. In turn, a sequence of random variables 
is uniformly integrable if it is bounded in $L^q$ for some $q>1$. Let us choose a $q>1$ such that $qp<\beta$. By Minkowski's inequality
we conclude that 
\begin{align*}
\E [| n^{d(Hp-1)} V_n(p) |^q] &\le\Bigg( n^{-d} \sum_{\boldsymbol{0} \le \boldsymbol{i} < \boldsymbol{n}} \Big(\E [ |n^{dH} \Delta_{1/n} X ( \boldsymbol{i}/n ) |^{qp} ] \Big)^{\frac{1}{q}} \Bigg)^q = \E [|n^{dH} \Delta_{1/n} X(\boldsymbol{0})|^{qp} ] \\
&= \E [|L([0,1]^d)|^{qp}] \Big(\int_{\R^d} |n^{d \alpha} \Delta_{1/n} g(\boldsymbol{s}/n)|^\beta \d \boldsymbol{s}\Big)^{qp/\beta} = O(1),
\end{align*}
where the last relation follows from \eqref{lim:Yn-Yinfty}.
Hence, the statistic in Theorem \ref{thm1}(ii) is uniformly integrable, and the proof is complete.

\subsection{Proof of Theorem~\ref{thm1}(iii)}

We start noticing that  under \textnormal{(H1)}, $g$ has continuous partial derivatives up to $d$-order in all 
 $\boldsymbol{s}=(s_1,\dots,s_d)\in \R^d$ with $s_i\neq 0$ for all $i=1,\dots,d$. 
Furthermore, 
\begin{equation}\label{lsjdfljsldh}
 (a): \int_{B_\rho (\boldsymbol{0})} | \partial^d g(\boldsymbol{s}) |^{\beta} \d \boldsymbol{s}<\infty \qquad \text{and}\qquad 
  (b): \int_{B_\rho (\boldsymbol{0})} | \partial^d g(\boldsymbol{s}) |^{p} \d \boldsymbol{s}<\infty,
\end{equation}
which follows from the estimate  $|\partial^d g (\boldsymbol{s})| \le C \| \boldsymbol{s} \|^{d(\alpha-1)}$ for all $\boldsymbol{s} \in B_\rho ( \boldsymbol{0} )$, and the fact that  $\int_{B_\rho (\boldsymbol{0})} \| \boldsymbol{s} \|^{dr(\alpha-1)} \d \boldsymbol{s} < \infty$ if and only if $dr(\alpha-1)+d-1>-1$. The latter condition is satisfied for $r=p$ and  $r=\beta$ since $1<\alpha + 1/\max(\beta,p)$. From \eqref{lsjdfljsldh}(b) and  $p\geq 1$, we deduce that  $\int_{B_\rho (\boldsymbol{0})} | \partial^d g(\boldsymbol{s}) | \d \boldsymbol{s}<\infty$ from which we conclude that 
\begin{equation}\label{sdlfjsldhf}
g([\boldsymbol{s},\boldsymbol{t}])= \int_{[\boldsymbol{s},\boldsymbol{t}]} \partial^d g(\boldsymbol{u})\d \boldsymbol{u},\qquad \text{for all } \boldsymbol{s}\leq \boldsymbol{t},
\end{equation}
where the left-hand side of \eqref{sdlfjsldhf} denotes the 
increments of $g$ over $[\boldsymbol{s},\boldsymbol{t}]$ defined in  \eqref{def:incrX}.  We now define a process $Y=(Y(\boldsymbol{t}))_{\boldsymbol{t}\in [0,1]^d}$ by 
\begin{equation}\label{sdlfjsdsh}
Y(\boldsymbol{t})= \int_{\R^d} \partial^d g(\boldsymbol{t} - \boldsymbol{s}) L(\d \boldsymbol{s}).
\end{equation}
It follows from \cite[Theorem 2.7]{rajput1989}, that $Y(\boldsymbol{t})$ is well-defined if and only if 
\begin{equation}\label{def:VL}
\int_{\R^d} V \big(\partial^d g(\boldsymbol{s}) \big) \d \boldsymbol{s} < \infty, \qquad \text{where } V (x) := \int_0^\infty \min \big(|xy|^2, 1 \big) \nu (\d y) \text{ for } x \in \R.
\end{equation}
Recall the estimate \eqref{ineq:estV}, where we have $V(x) \le C |x|^\theta$ for $|x| < 1$, whereas $V(x) \le C |x|^\beta$ for $|x| \ge 1$. 
By assumption (H1), there exists a $\rho>0$ such that $\partial^d g$ is bounded on $B_\rho^c (\boldsymbol{0})$ and is in 
$L^\theta ( B_\rho^c (\boldsymbol{0}) )$, and  $\partial^d g\in L^\beta( B_\rho (\boldsymbol{0}) )$, cf.\ \eqref{lsjdfljsldh}(a), which shows \eqref{def:VL}.

Next we will show existence of a  measurable and separable modification of $Y$ with values in the extended reals $[-\infty,\infty]$, and to this aim we let
 $L^\Phi$ denote the Musielak--Orlicz space of all $h:\R^d\to \R$ with 
\begin{equation}
 \Phi(h) := \int_{\R^d} \Big( \int_0^\infty \big( | y h(\boldsymbol{s}) |^2 \wedge 1 \big) \nu (\d y) \Big) \d \boldsymbol{s}<\infty 
\end{equation}
equipped with the $F$-norm 
\begin{equation}
\| h\|_\Phi= \inf\{c>0: \Phi(h/c)\leq 1\}. 
\end{equation}
Then, $L^\Phi$ is a separable linear metric space, and hence the mapping 
$\boldsymbol{t}\mapsto f_{\boldsymbol{t}}:= \partial^d g(\boldsymbol{t}-\cdot)$ from $[0,1]^d$ into $L^\Phi$ is measurable 
if $\boldsymbol{t}\mapsto \| g -f_{\boldsymbol{t}} \|_\Phi$ is measurable for all $g\in L^\Phi$. However, the latter follows directly from the joint measurability of $(\boldsymbol{s}, \boldsymbol{t})\mapsto \partial^d g(\boldsymbol{t}-\boldsymbol{s})$. Since the mapping $h\in L^\Phi$ into 
$\int_{\R^d} h(\boldsymbol{s}) L(\d \boldsymbol{s})\in L^0$ is continuous, cf.\ Theorem~3.3 of \cite{rajput1989}, it follows that 
 the mapping $\boldsymbol{t}\in [0,1]^d$ into $Y(\boldsymbol{t})\in L^0$ is measurable, from which we conclude that 
 there exists a measurable and separable modification of   $(Y(\boldsymbol{t}))_{\boldsymbol{t}\in [0,1]^d}$, cf.\ Theorem~3 of \cite{Cohn}.
 In the following $(Y(\boldsymbol{t}))_{\boldsymbol{t}\in [0,1]^d}$ will always denote  such  measurable and separable modification.

\bigskip
\noindent
{\it Step 1.} We now consider the integrability of 
$Y=(Y(\boldsymbol{t}))_{\boldsymbol{t}\in [0,1]^d}$ with respect to $\boldsymbol{t}$.  
 It follows from \cite[Theorem~3.1(i)]{braverman1998} that  $Y$ has  
sample paths in $L^p ([0,1]^d, \lambda^d)$ almost surely if  the following conditions hold:
\begin{equation}\label{cond1}
\| \partial^d g(\cdot - \boldsymbol{s}) \|_p := \Big( \int_{[0,1]^d} |\partial^d g (\boldsymbol{t}-\boldsymbol{s})|^p \d \boldsymbol{t} \Big)^{1/p} < \infty \qquad \text{for } \lambda^d \text{-almost every } \boldsymbol{s} \in \R^d;
\end{equation}
for some $c>0$ and $\delta'>0$,
\begin{align}\label{cond2}
\int_{\R^d} \nu \Big( \Big( \frac{c}{\| \partial^d g (\cdot - \boldsymbol{s}) \|_p} , \infty \Big) \Big)  \d \boldsymbol{s} < \infty \qquad \text{and} \qquad \int_{[0,1]^d} \sigma^p (\boldsymbol{t}) \d \boldsymbol{t} < \infty,
\end{align}
where 
\begin{equation*}
\sigma (\boldsymbol{t}) := \inf \{ \theta > 0 : \Phi ( \partial^d g (\boldsymbol{t} - \cdot)/\theta ) \le \delta' \},
\end{equation*}
and
\begin{equation}\label{cond3}
\int_{[0,1]^d} \Big( \int_{\R^d} \Big( \int_{c\sigma(\boldsymbol{t})/|\partial^d g(\boldsymbol{t}-\boldsymbol{s})|}^{c/\| \partial^d g(\cdot - \boldsymbol{s}) \|_p} |y \partial^d g(\boldsymbol{t}-\boldsymbol{s})|^p \nu (\d y) \Big) \d \boldsymbol{s} \Big) \d \boldsymbol{t} < \infty,
\end{equation}
where the inner integral in the last formula is set to be zero, if its lower limit of integration exceeds the upper limit.

The condition \eqref{cond1} holds because
$\partial^d g$ is bounded on $B_\rho^c (\boldsymbol{0})$ and $|\partial^d g(\boldsymbol{t})| \le C \| \boldsymbol{t} \|^{d(\alpha-1)}$ for all $\boldsymbol{t} \in B_\rho (\boldsymbol{0})$, where $\int_{B_\rho (\boldsymbol{0})} \| \boldsymbol{t} \|^{pd(\alpha-1)} \d \boldsymbol{s} < \infty$ if and only if $pd(\alpha-1) + d-1 > -1$. Next, let us verify the first condition in \eqref{cond2}. Let $\rho > 0$ be large enough. For $\boldsymbol{s} \in B_{2\rho}(\boldsymbol{0})$, use $\| \partial^d g(\cdot - \boldsymbol{s})\|_p \le C$,  furthermore, $\nu((1/C,\infty)) < \infty$.
For $\boldsymbol{s} \in B^c_{2\rho} (\boldsymbol{0})$, $\boldsymbol{t} \in [0,1]^d$, note that $| \partial^d g(\boldsymbol{t}-\boldsymbol{s}) | \le |\partial^d g(\boldsymbol{s}/2)|$, which leads to $\| \partial^d g(\cdot - \boldsymbol{s})\|_p \le |\partial^d g(\boldsymbol{s}/2)|$.
Finally, use that $\partial^d g \in L^\theta (B^c_{\rho}(\boldsymbol{0}))$ is bounded and $\nu ((y,\infty) )\le C y^{-\theta}$ for $y \ge 1$ to see that
$$
\int_{B^c_{2\rho}(\boldsymbol{0})} \nu \Big( \Big( \frac{c}{\| \partial^d g (\cdot - \boldsymbol{s}) \|_p} , \infty \Big) \Big)  \d \boldsymbol{s} \le  C \int_{B^c_{2\rho}(\boldsymbol{0})} |\partial^d g(\boldsymbol{s}/2)|^\theta \d \boldsymbol{s} < \infty.
$$
Note that $\Phi(\partial^d g(\boldsymbol{t}-\cdot))$, and hence $\sigma(\boldsymbol{t})$, both do not depend on $\boldsymbol{t} \in [0,1]^d$. With $V (x)$ as given in \eqref{def:VL}, we have that 
$$
\Phi(\partial^d g(\boldsymbol{t}-\cdot)) = \int_{\R^d} V (\partial^d g(\boldsymbol{s})) \d \boldsymbol{s} < \infty
$$
since $\alpha + 1/\beta > 1$. Hence, we conclude that the second condition in \eqref{cond2} holds.

Finally, we show \eqref{cond3}. Recall that $\rho$ is large enough so  that we have $|\partial^d g(\boldsymbol{t}-\boldsymbol{s})| \le |\partial^d g(\boldsymbol{s}/2)| \le C$ for $\boldsymbol{s} \in B^c_{2\rho}(\boldsymbol{0})$, $\boldsymbol{t} \in [0,1]^d$. We obtain
\begin{align*}
&\int_{[0,1]^d} \Big( \int_{B_{2\rho}^c (\boldsymbol{0})} \Big( \int_{c\sigma(\boldsymbol{t})/|\partial^d g(\boldsymbol{t}-\boldsymbol{s})|}^{c/\| \partial^d g(\cdot - \boldsymbol{s}) \|_p} |y \partial^d g(\boldsymbol{t}-\boldsymbol{s})|^p \nu (\d y) \Big) \d \boldsymbol{s} \Big) \d \boldsymbol{t} \\
&\quad \le  \int_{B_{2\rho}^c (\boldsymbol{0})} \Big( \frac{c^p}{\| \partial^d g(\cdot - \boldsymbol{s}) \|_p^p} \int_{[0,1]^d}  | \partial^d g(\boldsymbol{t}-\boldsymbol{s})|^p  \d \boldsymbol{t} \Big) \Big( \int_{C/|\partial^d g(\boldsymbol{s}/2)|}^{\infty}  \nu (\d y) \Big) \d \boldsymbol{s} \\
&\quad \le C \int_{B^c_{2\rho}(\boldsymbol{0})} |\partial^d g(\boldsymbol{s}/2)|^\theta \d \boldsymbol{s} < \infty.
\end{align*} 
We next deal with
\begin{align*}
I := \int_{[0,1]^d} \Big( \int_{B_{2\rho} (\boldsymbol{0})} \Big( \int_{c\sigma(\boldsymbol{t})/|\partial^d g(\boldsymbol{t}-\boldsymbol{s})|}^{c/\| \partial^d g(\cdot - \boldsymbol{s}) \|_p} |y \partial^d g(\boldsymbol{t}-\boldsymbol{s})|^p \nu (\d y) \Big) \d \boldsymbol{s} \Big) \d \boldsymbol{t}.
\end{align*}
If $p>\beta$, then for $\boldsymbol{t} \in [0,1]^d$, $\boldsymbol{s} \in B_{2\rho}(\boldsymbol{0})$,
\begin{align*}
&\int_0^\infty \1 \Big( \frac{c \sigma (\boldsymbol{t})}{|\partial^d g(\boldsymbol{t}-\boldsymbol{s})|} < y < \frac{c}{\| \partial^d g(\cdot - \boldsymbol{s}) \|_p} \Big) y^p \nu (\d y)\\
&\quad \le \int_0^1 y^p \nu (\d y) + \frac{c^p}{\| \partial^d g(\cdot - \boldsymbol{s}) \|_p^p} \int_1^\infty \nu (\d y) \le C \Big(1+ \frac{1}{\| \partial^d g(\cdot - \boldsymbol{s}) \|_p^p} \Big)
\end{align*}
and so 
\begin{align*}
I &\le C \int_{[0,1]^d} \Big( \int_{B_{2\rho}(\boldsymbol{0})} \Big( 1 + \frac{1}{\| \partial^d g(\cdot - \boldsymbol{s}) \|^p_p} \Big) |\partial^d g(\boldsymbol{t}- \boldsymbol{s})|^p \d \boldsymbol{s} \Big) \d \boldsymbol{t}\\
&= C \int_{B_{2\rho}(\boldsymbol{0})} ( \| \partial^d g(\cdot- \boldsymbol{s}) \|^p_p + 1 ) \d \boldsymbol{s}  < \infty.
\end{align*}
If $p \le \beta < \beta'$ with $\alpha+1/\beta'>1$, then for $\boldsymbol{t} \in [0,1]^d$, $\boldsymbol{s} \in B_{2\rho}(\boldsymbol{0})$,
\begin{align*}
&\int_0^\infty \1 \Big( \frac{c \sigma (\boldsymbol{t})}{|\partial^d g(\boldsymbol{t}-\boldsymbol{s})|} < y < \frac{c}{\| \partial^d g(\cdot - \boldsymbol{s}) \|_p} \Big) y^p \nu (\d y)\\
&\quad \le \int_0^1 \1 \Big( \frac{c \sigma (\boldsymbol{t})}{|\partial^d g(\boldsymbol{t}-\boldsymbol{s})|} < y < \frac{c}{\| \partial^d g(\cdot - \boldsymbol{s}) \|_p} \Big) y^{(p-\beta')+\beta'} \nu (\d y) + \frac{c^p}{\| \partial^d g(\cdot - \boldsymbol{s}) \|_p^p} \int_1^\infty \nu (\d y)\\ 
&\quad \le \Big( \frac{c\sigma(\boldsymbol{t})}{|\partial^d g(\boldsymbol{t}-\boldsymbol{s})|} \Big)^{p-\beta'} \int_0^1 y^{\beta'} \nu (\d y) + \frac{c^p}{\| \partial^d g(\cdot - \boldsymbol{s}) \|_p^p} \int_1^\infty \nu (\d y)\\
&\quad \le C \Big( |\partial^d g(\boldsymbol{t}-\boldsymbol{s})|^{\beta'-p} + \frac{1}{\| \partial^d g(\cdot - \boldsymbol{s})\|^p_p} \Big)
\end{align*}
and so 
\begin{align*}
I &\le C \int_{[0,1]^d} \Big( \int_{B_{2\rho}(\boldsymbol{0})} \Big( |\partial^d g(\boldsymbol{t}-\boldsymbol{s})|^{\beta'-p} + \frac{1}{\| \partial^d g(\cdot - \boldsymbol{s})\|^p_p} \Big) |\partial^d g(\boldsymbol{t}-\boldsymbol{s})|^p \d \boldsymbol{s} \Big) \d \boldsymbol{t}\\
&= C \int_{B_{2\rho}(\boldsymbol{0})} \Big( \int_{[0,1]^d} | \partial^d g(\boldsymbol{t} - \boldsymbol{s}) |^{\beta'} \d \boldsymbol{t} + 1 \Big) \d \boldsymbol{s} < \infty.
\end{align*}
We conclude that \eqref{cond3} holds. 
\bigskip

\noindent
{\it Step 2.}
In the following we will show that  for all $\boldsymbol{t}\in [0,1]^d$ we have almost surely 
\begin{equation}\label{sdlfjsdljf0}
 X([\boldsymbol{0},\boldsymbol{t}]) = \int_{[\boldsymbol{0},\boldsymbol{t}]} Y(\boldsymbol{u}) \d \boldsymbol{u}.
\end{equation}
	Note that the right-hand side of \eqref{sdlfjsdljf0}  is well-defined since $Y$ has sample paths in $L^p ( [0,1]^d, \lambda^d )\subseteq  L^1 ( [0,1]^d, \lambda^d )$. Choose a probability measure 
	$\kappa$ on $\R^d\times \R$ equivalent to $\lambda^d\otimes \nu$  and let $\eta$ denote the density of $\kappa$ with respect to $\lambda^d\otimes \nu$. 	According to  Theorem~5.1 and Remark~5.2(a)  of  \cite{Ros2018}  we  may choose three sequences
	$(\epsilon_j)_{j\in \N}$,  $(\Gamma_j)_{j\in \N}$ and $(\boldsymbol{\xi}_j)_{j\in \N}$, where $\boldsymbol{\xi}_j=(\boldsymbol{\xi}^1_j,\xi_j^2)\in \R^d\times \R$, such that 
		\begin{align}\label{sdlfjsdlfjshgs}
	Y(\boldsymbol{t}) = {}& \sum_{j=1}^\infty \epsilon_j \partial^d g(\boldsymbol{t}- \boldsymbol{\xi}_j^1)\xi_j^2 \1 ( \eta(\boldsymbol{\xi}_j)\leq \Gamma_j^{-1} ), \\	X(\boldsymbol{t}) ={}&  \sum_{j=1}^\infty \epsilon_j  g(\boldsymbol{t},\boldsymbol{\xi}_j^1) \xi^2_j\1 ( \eta(\boldsymbol{\xi}_j)\leq \Gamma_j^{-1} )
	\end{align}
	almost surely for all $\boldsymbol{t}\in [0,1]^d$. 
	Moreover, $(\boldsymbol{\xi}_j)_{j\in \N}$ is  an i.i.d.\ sequence of $\R^d\times \R$-valued random vectors with the common distribution $\kappa$,  $(\Gamma_j)_{j\in \N}$ is a sequence of partial sums of i.i.d.\ standard exponential random variables, 
	and  $(\epsilon_j)_{j\in \N}$ denotes an i.i.d.\ sequence of symmetric Bernoulli random variables, that is,  $\P(\epsilon_j=1)=\P(\epsilon_j=-1)=1/2$ for all $j\in \N$.  In addition, the three sequences $(\boldsymbol{\xi}_j)_{j\in \N}$, $(\Gamma_j)_{j\in \N}$ and $(\epsilon_j)_{j\in \N}$  are independent. Since conditionally on $(\boldsymbol{\xi}_j,\Gamma_j)_{j\in \N}$, the summands in \eqref{sdlfjsdlfjshgs} are independent and symmetric random elements with values in $L^1([0,1]^d,\lambda^d)$ and furthermore $Y$ has paths in $L^1([0,1]^d,\lambda^d)$, it follows by the It\^o--Nisio theorem, see \cite[Theorem~2.1.1]{KwaWoy}, 
	that the series \eqref{sdlfjsdlfjshgs} convergence in $L^1([0,1]^d,\lambda^d)$ with probability one. 
	In particular, for all $\boldsymbol{t}\in [0,1]^d$ we have with probability one 
	\begin{align}
\int_{[\boldsymbol{0},\boldsymbol{t}]} Y (\boldsymbol{u}) \d \boldsymbol{u} = {}& \sum_{j=1}^\infty 
\epsilon_j \Big(\int_{[\boldsymbol{0},\boldsymbol{t}]} \partial^d g(\boldsymbol{u}-\boldsymbol{\xi}_j^1) \d \boldsymbol{u}\Big) \xi^2_j \1 (\eta(\boldsymbol{\xi}_j)\leq \Gamma_j^{-1} )
\\ = {}& \sum_{j=1}^\infty 
\epsilon_j g([-\boldsymbol{\xi}_j^1,\boldsymbol{t}-\boldsymbol{\xi}_j^1]) \xi^2_j \1 (\eta(\boldsymbol{\xi}_j)\leq \Gamma_j^{-1} ) =X([\boldsymbol{0},\boldsymbol{t}]),
	\end{align}
	where the second equality follows by \eqref{sdlfjsldhf}. 
Hence the proof of \eqref{sdlfjsdljf0} is complete.
 \bigskip

\noindent
{\it Step 3.} For a given $p \ge 1$, we denote by $AC^p ([0,1]^d)$ the space of 
functions $\xi:[0,1]^d\to \R$ such that there is a function $\partial^d \xi \in L^p ([0,1]^d, \lambda^d)$ with 
\begin{equation}
 \xi ([{\bf 0},\boldsymbol{t}])
= \int_{[\boldsymbol{0}, \boldsymbol{t}]} \partial^d \xi (\boldsymbol{u}) \d \boldsymbol{u}, \qquad \text{for all } \boldsymbol{t} \in [0,1]^d.
\end{equation}
For $\xi \in AC^p([0,1]^d)$  let us prove that as $n \to \infty$,
\begin{equation}\label{lim:Vpxi}
n^{d(p-1)} V_n^\xi (p) := n^{d(p-1)} \sum_{\boldsymbol{0} \le \boldsymbol{i} < \boldsymbol{n}} | \Delta_{1/n} \xi (\boldsymbol{i}/n) |^p \to \int_{[0,1]^d} |\partial^d \xi (\boldsymbol{t})|^p \d \boldsymbol{t}.
\end{equation}
Firstly, assume that $\xi : \R^d \to \R$ has continuous partial derivatives up to the $(2d)$-th order at every point $\boldsymbol{t} \in \R^d$. We have that $n^d \Delta_{1/n} \xi (\boldsymbol{i}/n) = \partial^d \xi (\boldsymbol{i}/n) + r_n (\boldsymbol{i}/n)$, 
where $| r_n (\boldsymbol{i}/n) | \le C /n$ uniformly for all $\boldsymbol{0} \le \boldsymbol{i} < \boldsymbol{n}$. By Minkowski's inequality, 
\begin{align*}
\Big| \Big( n^{d(p-1)} V_n^\xi (p) \Big)^{1/p} - \Big( n^{-d} \sum_{\boldsymbol{0} \le \boldsymbol{i} < \boldsymbol{n}} | \partial^d \xi (\boldsymbol{i}/n) |^p \Big)^{1/p} \Big| \le \Big( n^{-d} \sum_{\boldsymbol{0} \le \boldsymbol{i} < \boldsymbol{n}} | r_n (\boldsymbol{i}/n) |^p \Big)^{1/p} = o (1)
\end{align*}
as $n \to \infty$.
By continuity of $\partial^d \xi$, we have that 
$$
n^{-d} \sum_{\boldsymbol{0} \le \boldsymbol{i} < \boldsymbol{n}} | \partial^d \xi (\boldsymbol{i}/n) |^p \to \int_{[0,1]^d} | \partial^d \xi (\boldsymbol{t}) |^p \d \boldsymbol{t} \qquad \text{as } n \to \infty.
$$
This proves \eqref{lim:Vpxi}.
Then, for general $\xi \in AC^p ([0,1]^d)$, $p\ge 1$, we  approximate $V_n^\xi (p)$ by $V_n^{\xi_m} (p)$, where $( \xi_m )$ is a sequence of functions having continuous partial derivatives up to the $(2d)$-th order at every point in $\R^d$. Indeed, the existence of such a sequence follows since continuous functions are dense in
$L^p([0,1]^d, \lambda^d)$.
 A combination of  \eqref{sdlfjsdljf0} and \eqref{lim:Vpxi} finishes the proof Theorem~\ref{thm1}(iii).

\subsection{Proof of Theorem \ref{thm2}(i)}

We denote by $Z$ the limiting variable
$$
\int_{[0,1]^k \times \R^{d-k} \times \R_0} \prod_{j=1}^k H_j (u_j) \prod_{j=k+1}^d \| g'_j (\cdot-x_j)\|^p_p |y|^p \Lambda^\ddagger (\d \boldsymbol{u}, \d \boldsymbol{x}, \d y)
$$
with $H_j (u) := \sum_{l \in \Z} |\Delta_1 h_j (l-u)|^p$, $u\in(0,1)$, $j=1,\dots,k$, and $\| g'_j(\cdot-x) \|_p := (\int_0^1 |g'_j(t-x)|^p \d t)^{1/p}$, $x \in \R$, $j=k+1,\dots,d$, where $\Lambda^\dagger$ is a Poisson random measure with intensity measure $\lambda^k \otimes \lambda^{d-k} \otimes \nu$ on $[0,1]^k \times \R^{d-k} \times \R_0$ introduced in Definition~\ref{sldjfsghsklhs}.  Then $|Z|<\infty$ almost surely if 
\begin{align*}
\int_{[0,1]^k \times \R^{d-k} \times \R_0} \min \Big( 1, \prod_{j=1}^k H_j (u) \prod_{j=k+1}^d \| g'_j (\cdot - x_j) \|_p^p |y|^p \Big) \d \boldsymbol{u} \d \boldsymbol{x} \nu (\d y) < \infty. 
\end{align*}
As in Theorem \ref{thm1}(i) for $d=1$, we have $\int_0^1 H_j(u) \d u< \infty$ since $\alpha_j+1/p\in(0,1)$, $j =1,\dots,k$. Hence, we only need to show that
\begin{equation}\label{exist5}
\int_{\R^{d-k} \times (0,\infty)} \min \Big( 1,\prod_{j=k+1}^d \| g'_j(\cdot - x_j)\|_p^p |y|^p \Big) \d \boldsymbol{x} \nu (\d y) < \infty.
\end{equation}
If $p \neq \theta < 2$, 
then for $|x| \le 1$, 
\begin{align*}
\int_0^\infty \min (1,|xy|^p) \nu (\d y) \le |x|^p \int_0^1 y^p \nu (\d y) + C \Big( |x|^p \int_1^{1/|x|} y^{p-\theta-1} \d y + \int_{1/|x|}^\infty y^{-\theta-1} \d y \Big) \le C |x|^{\min(p,\theta)},
\end{align*}
since $p>\beta$. On the other hand, if $\theta= 2$ then for $x \in \R$,
\begin{align*}
\int_0^\infty \min(1,|xy|^p) \nu (\d y) \le \int_0^\infty \min(1,|xy|^{\min(p,2)}) \nu (\d y) \le C |x|^{\min(p,2)},
\end{align*}
since $\min(p,2)> \beta$ and $\int_0^\infty y^2 \nu (\d y)< \infty$. This proves \eqref{exist5}, because for $\alpha_j +1/p>1$ and $|g'_j(s)| \ge |g'_j(t)|$ if $\rho \le |s| \le |t|$, it holds that
$$
\| g'_j (\cdot - x) \|_p \le C \1 (|x| < 2\rho) + |g'_j (x/2) | \1 (|x| \ge 2\rho), \qquad x \in \R,
$$ 
as in the proof of Theorem \ref{thm1}(iii) with $d=1$ (see also the verification of \eqref{cond1} in the proof of Theorem \ref{thm2}(iii)), moreover, $g'_j \in L^q ( (-\rho,\rho)^c )$, $q= \min (p,\theta)$, $j=k+1,\dots,d$. \\

\noindent
{\it Step 1}. 
Let $\nu (\R_0)<\infty$. We choose a small $\epsilon \in (0,1)$ and a large $m\in \N$. The way how $m$ depends on $\epsilon$ will be specified later. Now we decompose every $\Delta_{1/n} X(\boldsymbol{i}/n) = \tilde M_{n,\epsilon}(\boldsymbol{i}) + \tilde R'_{n,\epsilon}(\boldsymbol{i}) + \tilde R_{n,\epsilon}(\boldsymbol{i})$ so that 
\begin{align}\label{def:Mhat}
\tilde M_{n,\epsilon}(\boldsymbol{i})  &=  \int_{[0,1]^k \times \R^{d-k}} \prod_{j=1}^k \1 
( | i_j/n-s_j |<\epsilon ) \prod_{j=k+1}^d \1 ( |s_j| \le m ) \Delta_{1/n} g ( \boldsymbol{i}/n- \boldsymbol{s} ) L(\d \boldsymbol{s}) ,\\ \label{def:Rprimehat}
\tilde R'_{n,\epsilon}(\boldsymbol{i}) &=  \int_{\R^k \setminus [0,1]^k \times \R^{d-k}} \prod_{j=1}^k \1 (
| i_j/n-s_j |<\epsilon ) \prod_{j=k+1}^d \1 ( |s_j| \le m ) \Delta_{1/n} g ( \boldsymbol{i}/n - \boldsymbol{s} ) L(\d \boldsymbol{s}),\\ \label{def:Rhat}
\tilde R_{n,\epsilon}(\boldsymbol{i}) &=  \int_{\R^d}  \Big( 1 - \prod_{j=1}^k \1 (
| i_j/n -s_j | < \epsilon ) \prod_{j=k+1}^d \1 (|s_j| \le m) \Big)  \Delta_{1/n} g ( \boldsymbol{i}/n - \boldsymbol{s} ) L(\d \boldsymbol{s}) .
\end{align}
First, we will prove the stable convergence for the power variation statistics built from  $\tilde M_{n,\epsilon}(\boldsymbol{i})$ instead of the original increments $\Delta_{1/n} X(\boldsymbol{i}/n)$.
For this purpose, we use $\Lambda$ associated to $L$ by \eqref{ljsdlfjsdlj} and on $[0,1]^k \times \R^{d-k} \times \R_0$ having the representation $\Lambda = \sum_{l=1}^\infty \delta_{(\boldsymbol{V}_l^1, \boldsymbol{V}_l^2, J_l)}$ with $(\boldsymbol{V}_l^1,\boldsymbol{V}_l^2)=\boldsymbol{V}_l =(V_{l,1},\dots,V_{l,d})$ given in \eqref{sdljfsldhfs}. 
Particularly, we express the terms $\tilde M_{n,\epsilon}(\boldsymbol{i})$ as integrals with respect to $\Lambda$ on $[0,1]^k \times \R^{d-k} \times \R_0$:
\begin{align*}
\tilde M_{n,\epsilon}(\boldsymbol{i}) &= \int_{[0,1]^k \times \R^{d-k} \times \R_0} \prod_{j=1}^k \1 ( | i_j/n -x_j |<\epsilon ) \prod_{j=k+1}^d \1 (|x_j| \le m) \Delta_{1/n} g ( \boldsymbol{i}/n-\boldsymbol{x} ) y \, \Lambda(\d \boldsymbol{x},\d y)\\
&=\sum_{l=1}^{\infty} \prod_{j=1}^k \1 ( | i_j/n - V_{l,j} |<\epsilon ) \prod_{j=k+1}^d \1 (| V_{l,j} | \le m) \Delta_{1/n} g ( \boldsymbol{i}/n-\boldsymbol{V}_l ) J_l.
\end{align*}
We will prove that as $n \to \infty$ on the event
\begin{align*}
\Omega_\epsilon := \big\{ \omega \in \Omega : {}& \| \boldsymbol{V}_{l_1}^1 (\omega) - \boldsymbol{V}_{l_2}^1 (\omega) \|_\infty >2\epsilon   \text{ and } \boldsymbol{V}^1_{l_1}, \boldsymbol{V}^1_{l_2} (\omega) \in [\epsilon,1-\epsilon]^k  \\ 
{}&\text{ with } J_{l_1}(\omega),J_{l_2}(\omega)\neq 0 \text{ for all } l_1,l_2=1,2, \dots \big\}
\end{align*}
it holds
\begin{equation}\label{lim:Mj}
n^{\alpha_1 p + p -1} \sum_{\boldsymbol{0} \le \boldsymbol{i} < \boldsymbol{n}} | \tilde M_{n,\epsilon}(\boldsymbol{i}) |^p \overset{{\cal F}\textnormal{-d}}{\to} \sum_{l=1}^\infty |J_l|^p \prod_{j=1}^k H_j (U_{l,j}) \prod_{j=k+1}^d \| g'_j (\cdot - V_{l,j}) \|_p^p \, \1 ( | V_{l,j} | \le m ) = \tilde Z,
\end{equation}
where $(\boldsymbol{U}_l)_{l \in \N}$ with $\boldsymbol{U}_l = (U_{l,1},\dots,U_{l,k})$ is a sequence of independent ${\cal U}([0,1]^k)$-distributed random variables, defined on the extension of the underlying probability space $(\Omega,{\cal F},\P)$ and independent of the $\sigma$-algebra ${\cal F}$. 
To prove \eqref{lim:Mj}, we observe that on $\Omega_\epsilon$ every $\tilde M_{n,\epsilon}(\boldsymbol{i})$ satisfies either $|\tilde M_{n,\epsilon}(\boldsymbol{i})| = 0$ or $|\tilde M_{n,\epsilon}(\boldsymbol{i})| = |J_l \Delta_{1/n} g( \boldsymbol{i}/n - \boldsymbol{V}_l)| \prod_{j=k+1}^d \1 ( |V_{l,j}| \le m)$ for some $l=1,2,\dots$
Hence, it holds that on $\Omega_\epsilon$,
$$
\sum_{\boldsymbol{0} \le \boldsymbol{i} < \boldsymbol{n}} |\tilde M_{n,\epsilon}(\boldsymbol{i}) |^p = \tilde{V}_{n,\epsilon},
$$ 
where 
$$
\tilde V_{n,\epsilon} := \sum_{l=1}^\infty |J_l|^p \prod_{j=1}^k \sum_{i \in B_{n\epsilon}(\{n V_{l,j}\})} | \Delta_{1/n} g_j ( (i-\{ n V_{l,j} \})/n ) |^p \prod_{j=k+1}^d \sum_{0 \le i < n} | \Delta_{1/n} g_j (  i/n- V_{l,j} ) |^p \1 ( | V_{l,j} | \le m ).
$$
Since $\Omega_\epsilon \in {\cal F}$ then \eqref{lim:Mj} on $\Omega_\epsilon$ follows if we prove that 
\begin{equation}\label{sldfjsdlhs}
n^{\sum_{j=1}^k \alpha_j p + (d-k)(p - 1)} \tilde V_{n,\epsilon} \overset{{\cal F}\textnormal{-d}}{\to} \tilde Z \qquad \text{as } n \to \infty.
\end{equation}
To prove \eqref{sldfjsdlhs} we use the following arguments. On the left hand side of \eqref{sldfjsdlhs} each summand indexed by $l$ is a product of independent factors. As for these factors, we have
\begin{align*}
&n^{-1} \sum_{0 \le i < n} | n \Delta_{1/n} g_j ( i/n - V_{l,j} ) |^p \1 ( | V_{l,j} | \le m )\\ 
&\qquad\overset{\P}{\to} \int_0^1 | g'_j ( t-V_{l,j} ) |^p \d t \, \1 ( | V_{l,j} | \le m ) = \| g'_j ( \cdot -V_{l,j} ) \|_p^p \, \1 ( | V_{l,j} | \le m )
\qquad \text{as } n \to \infty,
\end{align*}
using Lemma~4.4 of \cite{basse2017} and $\alpha_j + 1/p >1$ for $j=k+1,\dots,d$, and
\begin{align*}
n^{\alpha_j p} \sum_{i \in B_{n\epsilon}(\{n V_{l,j}\})} | \Delta_{1/n} g_j ( (i-\{ n V_{l,j} \}) / n ) |^p \overset{{\cal F}\textnormal{-d}}{\to} \sum_{i \in \Z} |\Delta_1 h_j (i-U_{l,j})|^p = H_j(U_{l,j}),
\end{align*}
using the proof of Theorem \ref{thm1}(i) and $\alpha_j+1/p <1$ for $j=1,\dots,k$.
At last we note that as $\epsilon \downarrow 0$ together with $m\to \infty$,
\begin{align*}
\tilde Z\1 (\Omega_\epsilon) &\overset{\P}{\to}  \sum_{l=1}^\infty |J_l|^p \prod_{j=1}^k H_j (U_{l,j}) \prod_{j=k+1}^d \| g'_j (\cdot - V_{l,j}) \|_p^p\\
&= \int_{[0,1]^k \times \R^{d-k} \times \R_0} |y|^p \prod_{j=1}^k H_j (u_j) \prod_{j=k+1}^d \| g'_j(\cdot - x_j)\|^p_p \Lambda^\ddagger (\d \boldsymbol{u}, \d \boldsymbol{x}, \d y) =Z.  
\end{align*}

In the sequel let $m \to \infty$ so that $\epsilon m^{d-k} \downarrow 0$ as $\epsilon \downarrow 0$. Let us prove that the terms $\tilde R'_{n,\epsilon}(\boldsymbol{i})$ in \eqref{def:Rprimehat} satisfy
\begin{equation}\label{lim:Rprime}
\lim_{\epsilon \downarrow 0} \limsup_{n \to \infty} \P \Big( n^{\sum_{j=1}^k \alpha_j p + (d-k)(p - 1)} \sum_{\boldsymbol{0} \le \boldsymbol{i} < \boldsymbol{n}} |\tilde R'_{n,\epsilon} (\boldsymbol{i})|^p > \delta \Big) = 0
\end{equation}
for all $\delta>0$. 
The proof runs similarly to that of \eqref{lim:Mprime}. Recall that $\Lambda$ is associated to $L$ by \eqref{ljsdlfjsdlj} and use its representation $\Lambda = \sum_{l=1}^\infty \delta_{(\boldsymbol{V}'^1_{l,j},\boldsymbol{V}'^2_{l,j}, J'_l)}$ on $[-\epsilon,1+\epsilon]^k \setminus [0,1]^k \times [-m,m]^{d-k} \times \R_0$, analogous to that in \eqref{sdljfsldhfs}. Let $(\boldsymbol{V}'^1_{l,j},\boldsymbol{V}'^2_{l,j}) = \boldsymbol{V}_{l,j} = (V_{l,1},\dots,V_{l,d})$. Then, for $\bar p = \max(p,1)$, it holds that 
\begin{equation}\label{ineq:GVprimeJ}
\Big( n^{\sum_{j=1}^k \alpha_j p + (d-k)(p-1)} \sum_{\boldsymbol{0} \le \boldsymbol{i} < \boldsymbol{n}} |\tilde R'_{n,\epsilon}(\boldsymbol{i})|^p \Big)^{1/\bar p} \le \sum_{l=1}^\infty ( G_{n,\epsilon}(\boldsymbol{V}'_{l,j})|J'_l|^p )^{1/\bar p},
\end{equation}
furthermore, on the right hand side of \eqref{ineq:GVprimeJ}  every summand satisfies
$$
G_{n,\epsilon}(\boldsymbol{V}'_{l,j}) := n^{\sum_{j=1}^k \alpha_j p + (d-k)(p-1)} \sum_{\boldsymbol{0} \le \boldsymbol{i} < \boldsymbol{n}} | \Delta_{1/n} g ( \boldsymbol{i}/n - \boldsymbol{V}'_l ) |^p \prod_{j=1}^k \1 ( | i_j/n - V'_{l,j} | < \epsilon ) = O_\P(1)
$$
as $n\to\infty$.
The last property follows from
\begin{align}\label{lim:intgk}
n^{\alpha_j p + 1} \int_{|x|<1} |\Delta_{1/n} g_j (x)|^p \d x = O (1)
\end{align}
for $\alpha_j + 1/p \in (0,1)$, $j=1,\dots,k$, see the proof of Theorem \ref{thm1}(i), and
\begin{align}\label{lim:intgd0}
n^p \int_{|x|<2/n} |\Delta_{1/n} g_j(x)|^p \d x \le C n^p \int_{|x|<3/n} |g_j(x)|^p \d x \le C n^p \int_{|x|<3/n} |x|^{\alpha_j p} \d x = C n^{p-(\alpha_j p +1)} = o(1)
\end{align}
combined with
\begin{align}\label{lim:intgd1}
n^p \int_{|x|\ge 2/n} |\Delta_{1/n} g_j (x)|^p \d x &\le 
C \int_{ 
	|x| < \rho} |x|^{(\alpha_j-1)p} \d x + \int_{|x|>\rho} |g'_j(x)|^p \d x < \infty
\end{align}
for $\alpha_j + 1/p > 1$, $j=k+1,\dots,d$. From this we conclude \eqref{lim:Rprime} since the number of summands on the right hand side of \eqref{ineq:GVprimeJ} has mean $\lambda^{k} ([-\epsilon,1+\epsilon]^k \setminus [0,1]^k) \lambda^{d-k}([-m,m]^{d-k}) \nu(\R_0) = O(\epsilon m^{d-k}) = o(1)$.

Finally, consider the terms $\tilde R_{n,\epsilon}(\boldsymbol{i})$ in \eqref{def:Rhat}.
Let us prove that 
\begin{equation}\label{lim:Rhat}
\lim_{\epsilon \downarrow 0} \limsup_{n\to \infty} \P \Big( n^{\sum_{j=1}^k \alpha_j p + (d-k)(p - 1)} \sum_{\boldsymbol{0} \le \boldsymbol{i} < \boldsymbol{n}} |\tilde R_{n,\epsilon}(\boldsymbol{i})|^p>\delta \Big) = 0
\end{equation}
for all $\delta>0$. Without loss of generality we discuss the case, where 
\begin{align*}
\tilde R_{n,\epsilon} (\boldsymbol{i}) &= \int_{\R^d} \prod_{0 < j \le J} \1 ( | i_j/n - s_j |<\epsilon ) \prod_{J<j \le k} \1 ( | i_j/n - s_j | \ge \epsilon) \\ 
&\qquad \qquad\times \prod_{k<j \le K} \1(|s_j| > m) \prod_{K < j \le d} \1 (|s_j| \le m) \Delta_{1/n} g ( \boldsymbol{i}/n - \boldsymbol{s} ) L (\d \boldsymbol{s})
\end{align*}
for $0\le J < k<K \le d$. We note that if $k=d$ then the index set $k<j \le d$ is empty, but there is at least one $j$ in $J< j \le k$, whereas if $k< d$ then $J< j \le k$ can be empty but in that case there is at least one index $j$ in the set $k<j \le K$. Now we define a bounded function $\psi_j \in L^{\min(p,\theta)} (\R)$ so that
\begin{equation}\label{ineq:psij}
n | \Delta_{1/n} g_j ( i/n - x )| \1 ( | i/n- x| \ge \epsilon ) \le \psi_j (x)
\end{equation}
for all $x \in \R$, $0 \le i<n$, $J<j\le k$, and then we define $\pi_{J,K}(\boldsymbol{x}) = \prod_{J<j\le k}\psi_j (x_j) \prod_{k < j \le K} |g'_j (x_j/2)|$. We use $\Lambda$ associated to $L$ by \eqref{ljsdlfjsdlj} and set $\Lambda^1(\cdot) = \Lambda (\cdot \cap \{(\boldsymbol{x},y) \in \R^d \times\R_0 : \pi_{J,K}(\boldsymbol{x})|y| >1\})$ and for every $B \in {\cal B}_b (\R^d)$ set
$$
L^1 (B) = \int_{B \times \R_0} y \Lambda^1 (\d\boldsymbol{x}, \d y) \quad \text{and} \quad L^0 (B) = L(B) - L^1 (B).
$$
Then $L^0$ and $L^1$ are independent infinitely divisible random measures such that for every $B \in {\cal B}_b (\R^2)$,
\begin{align*}
\E\Big[ \e^{\i t L^0 (B)} \Big]&= \exp \Big( \int_{B \times \R_0} ( \e^{\i t y} - 1 - \i t y \1 (|y| \le 1) ) \1 ( \pi_{J,K}(\boldsymbol{x})| y | \le 1 ) \d \boldsymbol{x} \nu (\d y) \Big), \\
\E \Big[ \e^{\i t L^1 (B)} \Big] &= \exp \Big( \int_{B \times \R_0} ( \e^{\i t y} - 1 - \i t y \1 (|y| \le 1) ) \1 ( \pi_{J,K}(\boldsymbol{x}) | y | > 1 ) \d \boldsymbol{x} \nu (\d y) \Big), \qquad t \in \R.
\end{align*}
We decompose $\tilde R_{n,\epsilon} (\boldsymbol{i}) = n^{K-J} (\tilde Q^0_{n,\epsilon}(\boldsymbol{i}) + \tilde Q^1_{n,\epsilon}(\boldsymbol{i}))$ with
\begin{align*}
\tilde Q^l_{n,\epsilon}(\boldsymbol{i}) &:= \int_{\R^d} \prod_{0 < j \le J} \1 ( | i_j/n - s_j |<\epsilon ) \prod_{J<j \le k} \1 ( | i_j/n - s_j | \ge \epsilon) \\ 
&\qquad \qquad\times \prod_{k<j \le K} \1(|s_j| > m) \prod_{K < j \le d} \1 (|s_j| \le m) n^{K-J} \Delta_{1/n} g ( \boldsymbol{i}/n - \boldsymbol{s} ) L^l (\d \boldsymbol{s}), \qquad l=0,1.
\end{align*}
We claim that for all $\delta>0$,
\begin{equation}\label{lim:PtildeQ0}
\limsup_{n \to \infty} \P \Big( n^{\sum_{0 <j\le k} \alpha_j p + (d-k)(p-1) - (K-J)p} \sum_{\boldsymbol{0} \le \boldsymbol{i} < \boldsymbol{n}} |\tilde Q^l_{n,\epsilon} (\boldsymbol{i}) |^p > \delta \Big), \qquad l=0,1,
\end{equation}
are $0$ if $J < k$, and tend to zero as $m\to \infty$ if $J=k$. For $l=0$ it follows once we show that
\begin{equation}\label{EtildeQ0}
\sup_{n \in \N, \, \boldsymbol{0} \le \boldsymbol{i} < \boldsymbol{n}} n^{\sum_{0 < j \le J} \alpha_j p + J + (d-K)p} \E [|\tilde Q^0_{n,\epsilon} (\boldsymbol{i})|^p]
\end{equation}
is bounded since $\alpha_j + 1/p -1<0$, $J< j \le k$, in case $J<k$, and \eqref{EtildeQ0} tends to zero as $m \to \infty$ in case $J=k$. For this purpose, by Theorem 3.3 in \cite{rajput1989}, we need to show that 
\begin{align*}
&\int_{\R^d} \int_{\R_0} \Big| n^{\sum_{0 < j \le J} \alpha_j + J/p + d-K} \prod_{0 < j \le J} \1 ( | i_j/n - x_j |<\epsilon ) \prod_{J<j \le k} \1 ( | i_j/n - x_j | \ge \epsilon) \\ 
&\qquad \qquad\times \prod_{k<j \le K} \1(|x_j| > m) \prod_{K < j \le d} \1 (|x_j| \le m) n^{K-J} \Delta_{1/n} g ( \boldsymbol{i}/n - \boldsymbol{x} ) y \Big|^q \1 (\pi_{J,K} (\boldsymbol{x})|y| \le 1) \nu (\d y) \d \boldsymbol{x},
\end{align*}
is finite for $q=p$ and in addition $q=2$ if $p>2$. Here we rewrite $\Delta_{1/n} g(\boldsymbol{i}/n - \boldsymbol{x}) = \prod_{0<j\le d} \Delta_{1/n} g_j(i_j/n - x_j)$. For $q=p$ and in addition $q=2$ if $p>2$, 
$$
\int_{|x|<1} |n^{\alpha_j + 1/p} \Delta_{1/n} g_j (x)|^q \d x = O(1), \qquad n \to \infty,
$$
since $\alpha_j +1/p \in (0,1)$, $0<j\le J$, whereas 
$$
\int_{\R} |n \Delta_{1/n} g_j (x) |^q \d x = O(1), \qquad n \to \infty,
$$
since $K<j\le d$, as shown in \eqref{lim:intgk}, \eqref{lim:intgd0}, \eqref{lim:intgd1}. Next we use the dominating function $\psi_j$ in \eqref{ineq:psij} for the remaining factors indexed by $J<j \le k$ and $n |\Delta_{1/n} g_j(i/n-x)| \le |g'_j (x/2)|$ for $|x| >m$, $0\le i<n$, $k<j\le K$. Note that the resulting function $\prod_{J<j \le k} \psi_j (x_j) \prod_{k < j \le K} |g'_j (x_j/2)|  = \pi_{J,K} (\boldsymbol{x})$ on $\R^{k-J} \times ( [-m,m]^c )^{K-k}$ is bounded and $\min(p,\theta)$-th power integrable, which proves our statement about \eqref{EtildeQ0} and hence \eqref{lim:PtildeQ0} because $\int_{\R_0} |x y|^p \1 (|xy| \le 1) \nu (\d y) \le C |x|^{\min(p,\theta)}$ for $|x| \le 1$, where $p \neq \theta$ if $\theta <2$. 

Next, we show \eqref{lim:PtildeQ0} for $l=1$. We use \eqref{ineq:psij}, where $J<j \le k$, and $n |\Delta_{1/n} g_j(i/n-x)| \le |g'_j (x/2)|$ for $|x| >m$, $0\le i<n$, $k<j\le K$ to see that
\begin{align*}
\tilde Q^1_{n,\epsilon}(\boldsymbol{i}) &\le \int_{\R^d \times \R_0} \Big| y \prod_{0 < j \le J} \Delta_{1/n} g_j (i_j/n - x_j) \1 ( | i_j/n - x_j |<\epsilon ) \prod_{J<j \le k} \psi_j (x_j) \\ 
&\qquad \qquad\times \prod_{k<j \le K} g'_j (x_j/2) \1(|x_j| > m) \prod_{K < j \le d} \Delta_{1/n} g_j ( i_j/n - x_j )\1 (|x_j| \le m) \Big| \Lambda^1 (\d \boldsymbol{x}, \d y).
\end{align*}
We denote the term above on the right hand side by $\tilde {\cal Q}^1_{n,\epsilon} (\boldsymbol{i})$, but note that it does not depend on $i_j$, $J < j \le K$. Furthermore, for $\bar p = \max (p,1)$,
\begin{align*}
&\Big(n^{\sum_{0 < j \le k} \alpha_j p + (d-k)(p-1) - (K-J)p }  \sum_{\boldsymbol{0} \le \boldsymbol{i} < \boldsymbol{n}} |\tilde {\cal Q}^1_{n,\epsilon} (\boldsymbol{i})|^p \Big)^{1/\bar p}\\ 
&\qquad \le \int_{\R^d \times \R_0} \Big( |y|^p \prod_{0<j\le J} n^{\alpha_j p} \sum_{0 \le i_j < n} |\Delta_{1/n} g_j (i_j/n-x_j)|^p \1 (|i_j/n - x_j|<\epsilon)\\ 
&\qquad \qquad \qquad \qquad \times \prod_{J<j \le k} n^{\alpha_j p +1-p} |\psi_j (x_j)|^p \prod_{k<j\le K} |g'_j (x_j/2)|^p \1 (|x_j|>m)\\
&\qquad \qquad \qquad \qquad \times 
\prod_{K<j \le d} n^{p-1} \sum_{0 \le i_j < n} |\Delta_{1/n} g_j (i_j/n - x_j) |^p \1 (|x_j| \le m) \Big)^{1/\bar p} \Lambda^1 (\d \boldsymbol{x}, \d y) =: \tilde {\cal Z}^1_{n,\epsilon},
\end{align*}
where $\tilde {\cal Z}^1_{n,\epsilon}$ is well defined as integral with respect to Poisson random measure $\Lambda^1$ having intensity measure $\1 (\pi_{J,K}(\boldsymbol{x})|y|>1) \d \boldsymbol{x} \nu (\d y)$ since $\int_{\R_0} \1 (|xy|>1) \nu (\d y) \le C |x|^\theta$ for $|x| \le 1$ and $\pi_{J,K}(\boldsymbol{x})$ on $\R^{k-J} \times ([-m,m]^c)^{K-k}$ is bounded and $\theta$-th power integrable. Finally, following Step 1 we can show that $\tilde {\cal Z}^1_{n,\epsilon} = o_\P (1)$ as $n \to \infty$ if $J<k$, since $\alpha_j p + 1-p <0$, $J<j \le k$,  whereas if $J=k$, then
$\tilde {\cal Z}^1_{n,\epsilon}$ converges weakly to the integral 
\begin{align*}
&\int_{[0,1]^k \times \R^{d-k} \times \R_0} \Big( |y|^p \prod_{0<j \le k} H_j(u_j) \prod_{k<j\le K} |g'_j(x_j/2)|^p \1 (|x_j|>m)\\ 
&\qquad\qquad \qquad \qquad \qquad \times  \prod_{K<j\le d} \| g'_j (\cdot - x_j) \|_p^p \1 (|x_j| \le m) \Big)^{1/\bar p} \Lambda^{1,\ddagger} (\d \boldsymbol{u},\d \boldsymbol{x}, \d y) = \tilde {\cal Z}^1
\end{align*}
with respect to the Poisson random measure $\Lambda^{1,\ddagger} (\cdot)= \Lambda^\ddagger ( \cdot \cap \{ (\boldsymbol{u},\boldsymbol{x},y) \in [0,1]^k \times \R^{d-k} \times \R_0 : \pi_{k,K}(\boldsymbol{x})|y|>1 \})$ with intensity measure $\1 (\pi_{k,K}(\boldsymbol{x})|y|>1)\d \boldsymbol{u} \d \boldsymbol{x} \nu(\d y)$ on $[0,1]^k \times \R^{d-k} \times \R_0$
as $n \to \infty$, which further converges in probability to $0$ as $m \to \infty$. This finishes the proof of \eqref{lim:PtildeQ0}, hence of \eqref{lim:Rhat}.
Theorem \ref{thm2}(i) is proved in case $\nu (\R_0)<\infty$.

\bigskip

\noindent
{\it Step 2.} Let $\nu (\R_0) = \infty$. 
We aim to show that as $n\to\infty$, 
\begin{align}\label{ljsdfljsdflh}
&n^{\sum_{j=1}^k \alpha_j p + (d-k)(p - 1)} V^X_n(p)\nn\\
&\qquad\overset{{\cal F}\textnormal{-d}}{\to} \int_{[0,1]^k \times \R^{d-k} \times \R_0} |y|^p \prod_{j=1}^k H_i(u_j) \prod_{j=k+1}^d \| g'_j (\cdot - x_j) \|_p^p \Lambda^\ddagger (\d \boldsymbol{u}, \d \boldsymbol{x}, \d y) =: Z,
\end{align}
where the notation $V_n^X (p)$ is used to stress that  $V_n(p)$ is calculated for process $X$. For some small $\epsilon>0$, we decompose $\Delta_{1/n} X (\boldsymbol{i}/n) = \Delta_{1/n} X^{\le \epsilon} (\boldsymbol{i}/n) + \Delta_{1/n} X^{>\epsilon} (\boldsymbol{i}/n)$ following Step 2 of the proof of Theorem \ref{thm1}(i).
Since $\nu( [-\epsilon,\epsilon]^c ) < \infty$, we have that as $n \to \infty$, 
\begin{align}\label{sldjfsldjflsh}
&n^{\sum_{j=1}^k \alpha_j p + (d-k)(p - 1)} V_n^{X^{> \epsilon}}(p)\nn\\ 
&\qquad\overset{{\cal F}\textnormal{-d}}{\to} \int_{[0,1]^k \times \R^{d-k} \times [-\epsilon,\epsilon]^c} |y|^p \prod_{j=1}^k H_i(u_j) \prod_{j=k+1}^d \| g'_j (\cdot - x_j) \|_p^p \Lambda^\ddagger (\d \boldsymbol{u}, \d \boldsymbol{x}, \d y) =: Z^{>\epsilon},
\end{align}
as shown in Step~1. 
Since  
$Z^{>\epsilon} \overset{\P}{\to} Z$ as $\epsilon \downarrow 0$,
\eqref{ljsdfljsdflh} follows from \eqref{sldjfsldjflsh} if we show that for all $\delta >0$,
\begin{equation}\label{lim:P12}
\lim_{\epsilon \downarrow 0} \limsup_{n \to \infty} \P \Big( n^{\sum_{j=1}^k \alpha_j p + (d-k)(p - 1)} V_n^{X^{\le \epsilon}}(p) > \delta \Big) = 0
\end{equation}
using  \eqref{ineq0} and \eqref{ineq1}.
Furthermore,  \eqref{lim:P12} follows  by Markov's inequality, if we  prove that 
$$
\lim_{\epsilon \downarrow 0} \limsup_{n \to \infty} n^{\sum_{j=1}^k \alpha_j p + (d-k)(p - 1)+d} \E [|\Delta_{1/n} X^{\le \epsilon} (\boldsymbol{0})|^p] = 0.
$$
For the latter it suffices to show the convergence
\begin{equation}\label{lim:P22}
\lim_{\epsilon \downarrow 0} \limsup_{n \to \infty} \int_{\R^d} \int_{0<|y|\le\epsilon} \phi_p ( n^{\sum_{j=1}^k \alpha_j + (d-k)(1 - 1/p)+d/p} \Delta_{1/n} g(\boldsymbol{x}) y ) \nu (\d y) \d \boldsymbol{x} = 0,
\end{equation}
where 
$\phi_p(y) := |y|^p \1 (|y|>1) + |y|^2 \1(|y| \le 1)$ satisfies $\phi_p (y) \le |y|^p + |y|^2 \1 (p>2)$ for $y \in \R$, cf.\ Theorem 3.3 in \cite{rajput1989}. 
For $q=p$ and in addition $q=2$ if $p>2$, we have
\begin{align*}
&\int_{\R} |n^{\alpha_j + 1/p} \Delta_{1/n} g_j(x)|^q \d x\\ 
&\qquad\le \int_{|x|<1} |n^{\alpha_j + 1/p} \Delta_{1/n} g_j(x)|^q \d x + \int_{|x|>1} |n^{\alpha_j + 1/p-1} \psi_j(x)|^q \d x = O(1), \qquad j=1,\dots,k,
\end{align*}
and 
\begin{align*}
I_{n,j} (q) &:= \int_{\R} |n \Delta_{1/n} g_j (x)|^q \d x = O(1),  \qquad j=k+1,\dots,d,
\end{align*}
as shown in Step 1.
Finally, similarly to \eqref{ineq:numom0}, we get 
$\int_0^\epsilon y^p \nu (\d y) = O(\epsilon^{p-\beta}) = o (1)$
as $\epsilon \downarrow 0$, since $p>\beta$.
This completes the proof of \eqref{lim:P22} and \eqref{lim:P12}, and therefore the proof of Theorem \ref{thm2}(i).

\subsection{Proof of Theorem \ref{thm2}(ii)}

Let us first verify that the limiting constant
$m(p) := \E [ |L([0,1]^d)|^p ] (\prod_{j=1}^k I_j \prod_{j=k+1}^d I'_j)^{p/\beta}$ is finite. Indeed, for $j=1,\dots,k$, we have $I_j := \int_{\R} |\Delta_1 h_i(s)|^\beta \d s < \infty$ since $\alpha_j+1/\beta \in (0,1)$ as in Theorem \ref{thm1}(ii) in case $d=1$, whereas $I'_j = \int_{\R} |g'_j(s)|^\beta \d s <\infty$ follows from $| g'_j(s) | \le C |s|^{\alpha_j-1}$, $|s| < \rho$, and $g'_j \in L^\beta ( (-\rho,\rho)^c )$ for  $1< \alpha_j+1/\beta$, $j=k+1,\dots,d$.

Let us now prove that the convergence stated in Theorem \ref{thm2}(ii) holds in probability. Note that working on the assumption (H2) increments of $X$ can be approximated coordinate-wise since those of its kernel $g(\boldsymbol{s}) = \prod_{j=1}^d g_j(s_j)$ can be factorized $g( [\boldsymbol{s}, \boldsymbol{t}] ) = \prod_{j =1}^d (g_j (t_j)- g_j(s_j))$ for all $\boldsymbol{s}<\boldsymbol{t}$ in $\R^d$. We define the first approximation $(Z_n(\boldsymbol{i}))_{\boldsymbol{i} \in \Z^d}$ by
$$
Z_n(\boldsymbol{i}) := \int_{\R^d} \prod_{j=1}^k n^{H_j} \Delta_{1/n} g_j (i_j/n-s_j) \prod_{j=k+1}^d g'_j (i_j/n-s_j) L(\d \boldsymbol{s}).
$$
Then the above-stated convergence in probability follows using \eqref{ineq0}, \eqref{ineq1} if we prove
\begin{equation}\label{LLNZn}
n^{-d} \sum_{\boldsymbol{0} \le \boldsymbol{i} < \boldsymbol{n}} |n^{d -k+ \sum_{j=1}^k H_j} \Delta_{1/n} X(\boldsymbol{i}/n) - Z_n(\boldsymbol{i})|^p \overset{\P}{\to} 0 \qquad \text{and} \qquad
n^{-d} \sum_{\boldsymbol{0} \le \boldsymbol{i} < \boldsymbol{n}} |Z_n(\boldsymbol{i})|^p \overset{\P}{\to} m (p).
\end{equation} 
By Markov's inequality we deal with the first sequence with mean
$$
\E [ |n^{d-k+\sum_{j=1}^k H_j} \Delta_{1/n} X(\boldsymbol{0}) - Z_n(\boldsymbol{0})|^p ] = C ( \prod_{j=1}^k I_{n,j} R'_n )^{p/\beta},
$$
where
$$
I_{n,j} = \int_{\R} | n^{H_j} \Delta_{1/n} g_j(s)|^\beta \d s = \int_{\R} | n^{\alpha_j} \Delta_{1/n} g_j(s/n)|^\beta \d s = O(1), \qquad j = 1, \dots, k,
$$ 
follows from \eqref{lim:Yn-Yinfty} for $d=1$ and it remains to show
\begin{equation}\label{lim:Rn}
R'_n = \int_{\R^{d-k}} | \prod_{j=k+1}^d n \Delta_{1/n} g_j (s_j) -  \prod_{j=k+1}^d g'_j (s_j) |^\beta \d s_{k+1} \dots \d s_d = o(1).
\end{equation}
We rewrite the above integrand using the identity $\prod_{j=k+1}^d a_j -  \prod_{j=k+1}^d b_j = \sum_{\# J \ge 1} \prod_{j \in J} (a_j - b_j) \prod_{j \in J^c} b_j$, $\boldsymbol{a}, \boldsymbol{b} \in \R^{d-k},$
where the sum $\sum_{\# J \ge 1}$ is taken over all subsets $J \subseteq \{ k+1, \dots, d\}$ of cardinality $\# J \ge 1$. We thus reduce our task in \eqref{lim:Rn} to proving 
$$
\int_{\R} | n \Delta_{1/n} g_j (s) - g'_j (s) |^\beta \d s = o(1), \qquad j = k+1,\dots,d.
$$
We note that $n\Delta_{1/n} g_j (s) \to g'_j(s)$ for almost every $s$. Moreover, 
$|n \Delta_{1/n} g_j (s)| = |\int_0^1 g'_j(s+u/n) \d u| \le |g'_j(s/2)|$
for $|s| \ge 2\rho$ and $| n \Delta_{1/n} g_j(s) | \le C |s|^{\alpha_j-1}$ for $2/n\le |s| < 2\rho$. Hence, 
$\int_{|s| \ge 2/n} |n \Delta_{1/n} g_j(s)-g'_j(s)|^\beta \d s = o(1)$
by the dominated convergence theorem, whereas
$\int_{|s| < 2/n} |n \Delta_{1/n} g_j (s)|^\beta \d s \le C n^\beta \int_0^{3/n} s^{\alpha_j \beta} \d s = o(1)$
since $1<\alpha_j+1/\beta$, $j=k+1,\dots,d$. 

Now, we prove that the second convergence in \eqref{LLNZn} holds in $L^1$. Since for every $(i_{k+1},\dots,i_d) \in \Z^{d-k}$,
$$
( Z_n (\boldsymbol{i}) )_{(i_1,\dots,i_k) \in \Z^k} \overset{\rm fdd}{=} ( Z_n (i_1,\dots,i_k,0,\dots,0) )_{(i_1, \dots,i_k) \in \Z^k},
$$
it follows from 
\begin{equation}\label{lim:delta2X}
n^{- k} \sum_{0 \le i_1, \dots,i_k < n} | Z_n(i_1,\dots,i_k,0,\dots,0) |^p \overset{L^1}{\to} m(p).
\end{equation}
To show that the convergence  \eqref{lim:delta2X} holds in probability, we use the same arguments as in the proof of Theorem \ref{thm1}(ii). 
Using the scaling property of the $\beta$-stable random measure, we have that $
(Z_n (i_1,\dots,i_k,0,\dots,0))_{(i_1,\dots,i_k) \in \Z^k} \overset{\rm fdd}{=} (Y_n (i_1,\dots,i_k))_{(i_1,\dots,i_k) \in \Z^k},
$ 
and so 
$$
\sum_{0\le i_1,\dots,i_k<n} |Z_n(i_1,\dots,i_k,0,\dots,0)|^p \overset{\rm d}{=} \sum_{0 \le i_1,\dots,i_k < n} |Y_n(i_1,\dots,i_k)|^p,
$$
where
$$
Y_n (i_1,\dots,i_k) := \int_{\R^k} \prod_{j=1}^k n^{\alpha_j} \Delta_{1/n} g_j((i_j-s_j)/n) \prod_{j=k+1}^d g'_j (s_j) L(\d s_1, \dots, \d s_d).
$$
Next, we approximate $( Y_n(i_1,\dots,i_k))_{(i_1,\dots,i_k) \in \Z^k}$ by $Y_\infty = ( Y_\infty(i_1,\dots,i_k))_ {(i_1,\dots,i_k) \in \Z^k}$, where
$$
Y_\infty (i_1,\dots,i_k) := \int_{\R^d} \prod_{j=1}^k \Delta_{1} h_j (i_j-s_j) \prod_{j=k+1}^d g'_j (s_j) L(\d s_1, \dots, \d s_d),
$$
more specifically, we have that
\begin{align*}
&\E [|Y_n(0, \dots, 0)-Y_\infty(0, \dots, 0)|^p]\\ 
\qquad&= C \Big( \int_{\R} | \prod_{j=1}^k n^{\alpha_j} \Delta_{1/n} g_j(s_j/n)- \prod_{j=1}^k \Delta_1 h_j (s_j)|^\beta \d s_1 \dots \d s_k \prod_{j=k+1}^d \int_{\R} |g'_j (s)|^\beta \d s \Big)^{p/\beta} = o(1)
\end{align*}
using similar arguments to those in the proof of \eqref{lim:Rn} and \eqref{lim:Yn-Yinfty} for $d=1$. Hence, it follows that 
$$
n^{-k} \sum_{0 \le i_1,\dots,i_k < n} |Y_n(i_1,\dots,i_k)-Y_\infty(i_1,\dots,i_k)|^p \overset{\P}{\to} 0.
$$
Since the process 
$Y_\infty$
is a symmetric $\beta$-stable mixed moving average, by \cite[Theorem 3]{surg1993}, it is mixing, and hence ergodic. According to Birkhoff's theorem (see \cite[Theorem 10.6]{kalle2002}),
$$
n^{-k} \sum_{0 \le i_1,\dots,i_k < n} |Y_\infty(i_1,\dots,i_k)|^p \overset{\P}{\to} \E[ |Y_\infty(0,\dots,0)|^p],
$$
where $\E[ |Y_\infty (0, \dots,0)|^p] = m(p)$.
By \eqref{ineq0}, \eqref{ineq1} the sequence in \eqref{lim:delta2X} converges in probability. The sequence converges in mean if and only if it converges in probability and is uniformly integrable. 
The latter follows, because for some $q>1$ such that $qp<\beta$, by Minkowski's inequality, 
\begin{align*}
&\E \Big[\Big| n^{- k} \sum_{0 \le i_1,\dots,i_k < n} | Z_n(i_1,\dots,i_k,0,\dots,0) |^p \Big|^q \Big]\\ 
\qquad&\le \Big( n^{-k} \sum_{0 \le i_1,\dots,i_k < n} ( \E | Z_n(i_1,\dots,i_k,0,\dots,0) |^{qp} )^{1/q} \Big)^{q}
= \E [| Z_n(\boldsymbol{0}) |]^{qp}
= O(1).
\end{align*}
Similarly, $\E [|n^{d-k+\sum_{j=1}^k H_j} \Delta_{1/n} X(\boldsymbol{0})|^{qp}] = O(1)$, which completes the proof of Theorem \ref{thm2}(ii).

\subsection{Proof of Theorem \ref{thm2}(iii)}

The proof is analogous to that of Theorem \ref{thm1}(iii). It follows from \cite[Theorem 2.7]{rajput1989}, that the  random field $Y := ( Y ( \boldsymbol{t} ))_{\boldsymbol{t} \in [0,1]^d}$ given in \eqref{ljsdlfjhsgdhls} is well-defined if and only if 
\begin{align}\label{cond2:iii}
\int_{\R^d} V ( \partial^d g(\boldsymbol{s}) ) \d \boldsymbol{s} < \infty, \qquad \partial^d g(\boldsymbol{s}) := \prod_{i=1}^d g'_i (s_i), \qquad \boldsymbol{s} \in \R^d,
\end{align}
where 
$$
V (x) := \int_0^\infty \min(|xy|^2,1) \nu (\d y) \le C (|x|^\theta \1(|x|<1) + |x|^{\max(\beta,p)} \1(|x|\ge1)  ), \qquad x \in \R,
$$
as shown in \eqref{ineq:estV}. So \eqref{cond2:iii} follows from $g'_i \in L^\theta (\R) \cap L^{\max(\beta,p)} (\R)$, $i=1,\dots,d$, in case $\theta <\max(\beta,p)$ and from $g'_i \in L^{\max(\beta,p)} (\R)$, $i=1,\dots,d$, in case $\theta \ge \max(\beta,p)$. Note that (H2) implies that every $g'_i \in L^{q'} ((-\rho,\rho)^c)$ with $q' \ge \min(\theta,\max(\beta,p))$ and $|g'_i(s)| \le C |s|^{\alpha_i-1}$ for $|s| < \rho$ with $\alpha_i - 1 
> -1/\max(\beta,p) \ge -1/\min(\theta, \max(\beta,p))$, $i=1,\dots,d$. By the same arguments as in the proof of 
Theorem~\ref{thm1}(iii) we may choose a measurable and separable modification of $Y$, which also will be denoted $Y$. 


According to \cite[Theorem 3.1(i)]{braverman1998},   $Y$ has sample paths in $L^p ([0,1]^d, \lambda^d)$ almost surely if the conditions \eqref{cond1}, \eqref{cond2}, \eqref{cond3} hold.
For all $\boldsymbol{s} \in \R^d$, we have that
$\| \partial^d g(\cdot - \boldsymbol{s}) \|_p = \prod_{i =1}^d \| g'_i ( \cdot - s_i ) \|_p$,
where for 
$s \in \R$,
\begin{align*}
\| g'_i ( \cdot - s) \|_p := \Big( \int_{[0,1]} |g'_i(t-s)|^p \d t \Big)^{1/p} \le C \1 ( |s| < 2\rho ) + |g'_i(s/2)| \1 ( |s| \ge 2\rho ) \le C,
\end{align*}
because $|g'_i(t)|\le C|t|^{\alpha_i-1}$ for $|s|<3\rho$ with $\alpha_i - 1>-1/p$ and $|g'_i (s)| \ge |g'(t)|$ for  $1<\rho \le |s| \le |t|$, $i=1,\dots,d$.  
We conclude that condition \eqref{cond1} holds.

Next, let us verify the first condition in \eqref{cond2}. From above it follows that
\begin{align}\label{cond21}
\int_{\R^d} \nu \Big( \Big( \frac{c}{\|\partial^d g(\cdot - \boldsymbol{s})\|_p},\infty \Big) \Big) \d \boldsymbol{s} &\le C \int_{\R^d} \| \partial^d g(\cdot - \boldsymbol{s})\|_p^\theta \d \boldsymbol{s}\\
&\le C \int_{\R^d} \prod_{i=1}^d (\1(|s_i| < 2\rho) + |g'_i(s_i/2)|^\theta \1 (|s_i|\ge 2\rho) ) \d \boldsymbol{s} < \infty\nn
\end{align}
since $g'_i \in L^\theta ((-\rho,\rho)^c)$, $i=1,\dots,d$.
Note that 
$$
\Phi(\partial^d g(\boldsymbol{t}-\cdot)) = \int_{\R^d} V (\partial^d g(\boldsymbol{s})) \d \boldsymbol{s} < \infty,
$$ 
see \eqref{cond2:iii}, hence both $\Phi (\partial^d g(\boldsymbol{t}-\cdot))$ and $\sigma(\boldsymbol{t})$ do not depend on $\boldsymbol{t} \in [0,1]^d$. We conclude that the second condition in \eqref{cond2} holds.

For $0<c_0<c_1$, decompose $\int_{c_0}^{c_1} y^p \nu (\d y) = I_0 + I_1$, where
\begin{align*}
I_1 &:= \int_1^\infty \1(c_0<y<c_1) y^p \nu (\d y) \le \int_1^\infty \1 (c_0<y<c_1) y^{p-\theta-1} \d y\\
&\le C(c_0^{p-q} \1(p<q) + c_1^{p-q}\1 (p>q) +  \1 (p=q) )
\end{align*}
with $q=\min(\theta,\max(\beta,p))$ in case $p\neq \theta$, $\theta<2$ and $I_1 \le C$ in case $p=\theta=2$ and
\begin{align*}
I_0 := \int_0^1 \1(c_0<y<c_1) y^p \nu (\d y) &\le C \int_0^1 \left( \1(\beta<p) + \1 (p \le \beta) \1 (c_0<y) \right) y^p \nu (\d y)\\ 
&\le C ( \1(\beta<p) + \1(p \le \beta) c_0^{p-\beta'} )
\end{align*}
with $\beta' > \beta$ chosen so that $\min(\alpha_1,\dots,\alpha_d) + 1/\beta' > 1$. 
Therefore, the last condition \eqref{cond3} follows from 
\begin{align*}
&\int_{\R^d} (\| \partial^d g(\cdot-\boldsymbol{s})\|_{\beta'}^{\beta'} \1 (p \le \beta) + \| \partial^d g(\cdot-\boldsymbol{s})\|_p^p \1 (\beta<p)\\ 
&\qquad+ \| \partial^d g(\cdot - \boldsymbol{s})\|_q^q \1 (p<q) + \| \partial^d g(\cdot - \boldsymbol{s})\|_p^q \1(p>q) + \| \partial^d g(\cdot - \boldsymbol{s})\|_p^p \1(p=q) ) \d \boldsymbol{s} < \infty.
\end{align*}
To end the proof recall that $g'_i \in L^{q'}((-\rho,\rho)^c)$ with $q' \ge q$ and $|g'_i(s)| \le C|s|^{\alpha_i-1}$ for $|s|<\rho$ with $\alpha_i-1>
-1/\max(\beta,p) \ge -1/q$ , $i=1,\dots,d$.

\section{Appendix}

Let us verify that imposed Assumptions ($g$), ($\theta$) and  ($\beta$) for some $0<\theta \le 2$, $0 \le \beta < 2$ ensure the existence of the random field $X$. 
From \cite[Theorem 2.7]{rajput1989} it follows that the stochastic integral for $\boldsymbol{t} \in \R^d$ on the r.h.s.\ of \eqref{def:X} exists if and only if
\begin{equation}\label{Xexist}
\int_{\R^d}V (g(\boldsymbol{t},\boldsymbol{u})) \d \boldsymbol{u} < \infty 
\qquad \text{with }
V(x) := \int_0^\infty \min ( |xy|^2, 1 ) \nu (\d y) \text{ for } x \in \R,
\end{equation}
when $\nu$ is  a symmetric L\'evy measure on $\R$. Let us first show that Assumptions ($\beta$) and  ($\theta$) imply the following important estimate: there is a constant $C>0$ such that
\begin{equation}\label{ineq:estV}
V(x) \le C ( |x|^\theta \1 (|x| \le 1) + |x|^\beta \1 (|x| > 1)  ).
\end{equation}
Set $\bar \nu(y) := \nu (\{u\in \R_0 : u \ge y \})$ for $y>0$. If $\theta < 2$, then 
$y^{\theta} \bar \nu (y) \le C$ for $y \ge 1$, that is
$
\int_{1}^\infty f(u) \nu (\d u) \le C \int_{1}^\infty f(u) u^{-\theta-1} \d u
$ 
with $f(u)=\1(u\ge y)$, $u\in \R$, for $y\ge 1$, and 
the inequality remains valid by monotone approximation
for $f : [1,\infty) \to [0,\infty)$ non-decreasing. 
Hence, 
\begin{align}
V(x) 
&\le C \Big( x^2 + \int_{1}^\infty \min ( |xy|^2,1 ) y^{-\theta-1} \d y \Big) \nn \\
&\le C \Big( x^2 + x^2 \int_{1}^{\frac{1}{|x|}} y^{1-\theta} \d y + \int_{\frac{1}{|x|}}^\infty y^{-\theta-1} \d y \Big) \le C |x|^\theta \label{ineq:VL}
\end{align}
for $|x| \le 1$ if $\theta < 2$, whereas $V(x) \le C |x|^2$ for $x \in \R$ if $\theta = 2$. 

Furthermore, if $\beta>0$, then $y^\beta \bar \nu (y) \le C$ for $0<y<1$. For $0<\epsilon_0<\epsilon_1<1$, 
\begin{align*}
\int_{\epsilon_0}^{\epsilon_1} y^2 \nu (\d y) = - \int_{\epsilon_0}^{\epsilon_1} u^2 \bar \nu (\d u) = \epsilon_0^2 \bar \nu (\epsilon_0) - \epsilon_1^2 \bar  \nu (\epsilon_1) + 2 \int_{\epsilon_0}^{\epsilon_1} u^{1 -\beta} u^\beta \bar \nu (u) \d u,
\end{align*}
and so as $\epsilon_0 \to 0$,
\begin{equation}\label{ineq:numom0}
\int_0^{\epsilon_1} y^2 \nu (\d y) \le C \epsilon_1^{2-\beta}.
\end{equation}
Hence, 
$$
V (x) \le C \Big( |x|^2 \int_0^{\frac{1}{|x|}} |y|^2 \nu (\d y) + \int_{\frac{1}{|x|}}^\infty \nu (\d y) \Big)  \le C |x|^\beta
$$
for $|x| > 1$  if $\beta>0$,  
whereas
$V(x) \le C$ for $x \in \R$ if $\beta = 0$. This completes the proof of \eqref{ineq:estV}, and if moreover Assumption ($g$) holds, that of \eqref{Xexist}. We conclude that $X$ is well-defined.

\begin{center}
\Large{Acknowledgments}
\end{center}
The authors are grateful to an anonymous referee for useful comments.
Vytaut\.e Pilipauskait\.e and Mark Podolskij  gratefully acknowledge financial support
from the project ``Ambit fields: Probabilistic properties and statistical inference'' funded by Villum Fonden. Also, Vytaut\.e Pilipauskait\.e and Mark Podolskij gratefully acknowledge financial support of ERC Consolidator Grant 815703 ``STAMFORD: Statistical Methods for High Dimensional Diffusions''.


%

\end{document}